\documentclass[11pt]{article}

\usepackage{color}
\usepackage{latexsym}
\usepackage{amssymb}
\usepackage{graphicx}
\usepackage{amsmath}

\newtheorem{Theorem}{Theorem}[part]
\newtheorem{Definition}{Definition}[part]
\newtheorem{Proposition}{Proposition}[part]

\newtheorem{Lemma}{Lemma}[part]
\newtheorem{Corollary}{Corollary}[part]
\newtheorem{Remark}{Remark}[part]

\def \tri{\blacktriangleright}

\def \Frac{\displaystyle\frac}

\def \N{\mathbb{N}}
\def \R{\mathbb{R}}

\def \E{\mathbb{E}}
\def \F{\mathbb{F}}

\def \P{\mathbb{P}}

\def \T{\mathbb{T}}

\def \Ac{{\cal A}}

\def \Dc{{\cal D}}

\def \Fc{{\cal F}}

\def \Hc{{\cal H}}
\def \Ic{{\cal I}}

\def \Lc{{\cal L}}

\def \Oc{{\cal O}}

\def \Tc{{\cal T}}

\def \eps{\varepsilon}

\def \ep{\hbox{ }\hfill$\Box$}

\def\Dt#1{\Frac{\partial #1}{\partial t}}
\def\Dtp#1{\Frac{\partial #1}{\partial t'}}

\def\Dx#1{\Frac{\partial #1}{\partial x}}
\def\Dxp#1{\Frac{\partial #1}{\partial x'}}

\def\reff#1{{\rm(\ref{#1})}}

\def\beqs{\begin{eqnarray*}}
\def\enqs{\end{eqnarray*}}
\def\beq{\begin{eqnarray}}
\def\enq{\end{eqnarray}}

\begin{document}

\title{Impulse control problem on finite horizon \\ with  execution delay}

\author{Benjamin BRUDER
            \\\small   Laboratoire de Probabilit{\'e}s et
             \\\small  Mod{\`e}les Al{\'e}atoires, CNRS, UMR 7599
             \\\small  Universit{\'e} Paris 7 Diderot
             \\\small  and SGAM, Soci\'et\'e G\'en\'erale
             \\\small  bruder@math.jussieu.fr
             \and
         Huy{\^e}n PHAM
             \\\small  Laboratoire de Probabilit{\'e}s et
             \\\small  Mod{\`e}les Al{\'e}atoires, CNRS, UMR 7599
             \\\small  Universit{\'e} Paris 7 Diderot
              \\\small  and Institut Universitaire de France
              \\\small  pham@math.jussieu.fr}

%\date{}

\maketitle

\begin{abstract}
We consider impulse control problems in finite horizon for diffusions with  decision lag and execution delay.  The new feature is that our general framework deals with the important  case  when several consecutive orders may be decided before the  effective execution of the first one.   
This is motivated by financial applications in the trading of illiquid assets such as hedge funds.  
We show that the  value functions for such control problems satisfy a suitable version of dynamic programming principle in finite dimension, which takes into account the past dependence of state process  through the pending orders. The corresponding  Bellman partial differential equations (PDE) system  is derived, and  exhibit some  peculiarities on the coupled equations, domains and boundary conditions.  We prove a unique characterization of the value functions to this nonstandard PDE system  by means of viscosity solutions. We then provide an algorithm to find  the value functions and the optimal control.  This easily implementable algorithm involves backward and forward iterations on the domains and the value functions, which appear in turn
as original arguments in the proofs for the  boundary conditions and uniqueness results.
\end{abstract}

\vspace{13mm}

\noindent {\bf Key words~:} Impulse control, execution delay, diffusion processes, dynamic programming, viscosity solutions, comparison principle.

\vspace{5mm}

\noindent {\bf MSC Classification (2000):} 93E20, 62L15, 49L20, 49L25.

\newpage

\tableofcontents

\newpage

\section{Introduction}

\setcounter{equation}{0} \setcounter{Assumption}{0}
\setcounter{Theorem}{0} \setcounter{Proposition}{0}
\setcounter{Corollary}{0} \setcounter{Lemma}{0}
\setcounter{Definition}{0} \setcounter{Remark}{0}

In this paper, we  consider a general impulse control problem  in finite horizon of a diffusion process $X$,  with intervention lag and execution delay.
This means that  we may intervene on the diffusion system at any times $\tau_i$
separated at least by some fixed positive lag $h$, by  giving some impulse $\xi_i$ based on the information at $\tau_i$.
However, the execution of the impulse decided at $\tau_i$  is carried out with delay $mh$, $m$ $\geq$ $1$, i.e. it is implemented at time $\tau_i+mh$, moving the system from $X_{(\tau_i+mh)^-}$ to $\Gamma(X_{(\tau_i+mh)^-},\xi_i)$.   The objective is to maximize over impulse controls $(\tau_i,\xi_i)_i$ the expected total profit on finite horizon $T$,  of  the form
\beqs
\E\Big[  \int_0^T f(X_t) dt  + g(X_T) +  \sum_{\tau_i+mh\leq T}  c(X_{(\tau_i+mh)^-},\xi_i) \Big].
\enqs
Such formulations  appear  naturally  in  decision-making problems in economics and finance.  In many situations,  firms or investors face regulatory
delays (delivery lag),
which may be significant, and thus need to be taken into account  when management strategies  are decided in an uncertain environment.
Problems where  firm's investment are subject to delivery lag can be found in the real options literature, for example in  \cite{barstr96} and \cite{alvkep02}.
In financial market  context,  execution delay is related to liquidity risk (see e.g. \cite{subjar01}),  and occurs with
transaction, which requires heavy preparatory work as for hedge funds. Indeed, hedge funds frequently hold illiquid assets, and need some time to find a counterpart to buy or sell them. Furthermore, this notice period gives the hedge fund manager a reasonable investement horizon.

From a mathematical viewpoint, it is well-known that   impulse control problems without delay, i.e. $m$ $=$ $0$,  lead to variational
partial differential equations (PDE), see e.g.
the books \cite{benlio82} and  \cite{okssul05}.  Impulse control problems in the presence of delay were studied in  \cite{rob76}
for $m$ $=$ $1$, that is when  no more than one pending order is allowed at any time.  In this case, it is shown that the delay problem may be transformed
into a no-delay impulse control problem.   The paper \cite{bayega06} also considers  the case $m$ $=$ $1$,  but when the value of the impulse is chosen at the time of execution, and on infinite horizon, and these two conditions are crucial in the proposed probabilistic resolution.   We mention also  the
works  \cite{barsul95} and recently \cite{okssul06}, which study impulse problems in infinite horizon with arbitrary number of pending orders,
but under restrictive assumptions on the controlled state process, like  (geometric) L\'evy  process for $X$  and (multiplicative) additive intervention
operator $\Gamma$.  In this case,  the problem is  reduced to a finite-dimensional one where the value functions with pending orders are directly related to the value function without order.

The main contribution of this paper  is  to provide a  theory of impulse control problems with delay  on finite horizon
in a fairly general diffusion framework that deals with the important case in applications when  the number of pending orders is finite, but not restricted to
one, i.e. $m$ $\geq$ $1$.   Our chief goal is  to obtain a unique  tractable PDE characterization of the value functions for such problems.  As usual in stochastic control problems, the first  step is  the derivation of a dynamic programming principle (DPP).  We show a suitable version of DPP, which takes into account the past dependence  of the controlled diffusion  via the  finite number of pending orders.  The corresponding  Bellman PDE system
reveals some nonstandard features  both on the form of the differential operators and their domains, and on the boundary conditions.
Following the modern approach to stochastic control, we prove  that the value functions are viscosity solutions to this Bellman PDE system, and we also
state comparison principles, which allows to obtain a unique PDE characterization.  From this PDE  representation,  we provide an easily
implemented algorithm  to compute the value functions, and so as byproducts the optimal impulse control.  This algorithm involves forward and backward iterations on  the value functions and on the domains, and appear actually as original  arguments in the proofs for the boundary conditions and comparison principles.

The rest of the paper is organized as follows.  In Section 2, we formulate the control  problem and introduce the  associated value functions. Section 3 deals
with  the dynamic programming principle in this general framework. We then state in Section 4 the unique PDE viscosity characterization for the value functions.  In Section 5, we provide an algorithm for computing the value functions and the optimal impulse control.  Finally, Section 6 is  devoted to the proofs of  results in this paper.

\section{Problem formulation}

\setcounter{equation}{0} \setcounter{Assumption}{0}
\setcounter{Theorem}{0} \setcounter{Proposition}{0}
\setcounter{Corollary}{0} \setcounter{Lemma}{0}
\setcounter{Definition}{0} \setcounter{Remark}{0}

\subsection{The control problem}

Let $(\Omega,\Fc,\P)$ be a complete probability space equipped with a filtration $\F$ $=$ $(\Fc_t)_{t\geq 0}$ satisfying the usual conditions, and $W$ $=$
$(W_t)_{t\geq 0}$ a  standard $n$-dimensional Brownian motion.

An impulse control is a double sequence $\alpha$ $=$ $(\tau_i,\zeta_i)_{i\geq 1}$, where  $(\tau_i)$ is an increasing sequence of $\F$-stopping times,    and $\xi_i$ are $\Fc_{\tau_i}$-measurable random variables valued in $E$.  We require that
$\tau_{i+1}-\tau_i$ $\geq$ $h$ a.s., where $h$ $>$ $0$ is a fixed time lag between two decision times, and we assume that $E$, the set of impulse values,
is a compact subset of $\R^q$.  We denote by $\Ac$ this set of impulse controls.

In absence of  impulse executions,  the system valued in $\R^d$ evolves according to~:
\beq \label{eqX0}
dX_s &=& b(X_s) ds + \sigma(X_s) dW_s,
\enq
where   $b$ $:$ $\R^d$ $\rightarrow$ $\R^d$ and $\sigma$ $\rightarrow$ $\R^{d\times n}$ are Borel functions on $\R^d$,
satisfying usual Lipschitz conditions.   The interventions are  decided at times $\tau_i$ with impulse values $\xi_i$ based on the information at these dates,
however they are executed with delay at times  $\tau_i+mh$, moving the system from $X_{(\tau_i+mh)^-}$ to
$X_{(\tau_i+mh)}$ $=$ $\Gamma(X_{(\tau_i+mh)^-},\xi_i)$. Here
$\Gamma$ is a mapping from $\R^d\times E$ into  $\R^d$, and we assume that  $\Gamma$ is continuous, and satisfies  the linear growth condition~:
 \beq \label{hypgam}
\sup_{(x,e)\in\R^d\times E} \frac{ |\Gamma(x,e)|}{1+|x|}  & < & \infty.
 \enq
Given an impulse control $\alpha$ $=$ $(\tau_i,\xi_i)_{i\geq 1}$ $\in$ $\Ac$,
and an initial condition $X_0$ $\in$ $\R^d$,  the controlled process $X^\alpha$
is then defined as the solution to the s.d.e.~:
 \beq \label{eqX}
 \hspace{-7mm} X_s  &=& X_0 + \int_0^s  b(X_u) du + \int_0^s \sigma(X_u) dW_u +
  \sum_{\tau_i+mh \leq s} \big( \Gamma(X_{(\tau_i+mh)^-},\xi_i) -  X_{(\tau_i+mh)^-}\big).
 \enq
We now fix a finite horizon $T$ $<$ $\infty$,  and in order to avoid trivialities,
 we assume $T-mh$ $\geq$ $0$.
Using standard arguments based on Burkholder-Davis-Gundy's inequality, Gronwall's lemma and \reff{hypgam}, we easily check that
\beq \label{estimX}
\E\big[ \sup_{s\leq T} |X_s^{\alpha}| \big] & < &  \infty.
\enq

Given an impulse control $\alpha$ $=$ $(\tau_i,\xi_i)_{i\geq 1}$ $\in$ $\Ac$,
we consider the total profit  at horizon $T$,  defined by~:
 \beqs
 \Pi(\alpha) &=&  \int_0^T f(X_s^{\alpha}) ds + g(X_T^{\alpha}) + \sum_{\tau_i + mh \leq T} c(X_{(\tau_i+mh)^-}^\alpha,\xi_i),
 \enqs
 and we assume that the running profit function $f$,   the terminal profit function $g$,  and the executed cost function
 $c$  are continuous, and satisfy the linear growth condition~:
 \beq \label{hypfgc}
 \sup_{(x,e)\in\R^d\times E} \frac{ |f(x)| + |g(x)| + |c(x,e)|}{1+|x|}  & < & \infty.
 \enq
This ensures with \reff{estimX} that  $\Pi(\alpha)$ is integrable, and we can  define  the control problem~:
\beq \label{defcon}
V_0 &=& \sup_{\alpha\in\Ac} \E\big[ \Pi(\alpha) \big].
\enq

\noindent {\bf Financial example}

\noindent  Consider a financial market consisting of a money market account
yielding a constant interest rate $r$, and a risky asset (stock) of price process $(S_t)_t$ governed by~:
\beqs
dS_t &=& \beta(S_t) dt + \gamma(S_t) dW_t.
\enqs
We denote by $Y_t$ the number of shares
in the  stock, and by $Z_t$ the amount of money (cash holdings)  held by the investor at time $t$. We assume that the investor can only trade discretely, and her orders are executed with delay. This is modelled through an impulse control $\alpha$ $=$ $(\tau_i,\xi_i)_{i\geq 1}$ $\in$ $\Ac$, where $\tau_i$ are the decision times, and $\xi_i$ are the numbers of stock purchased if $\xi_i$ $\geq$ $0$ or selled if $\xi_i$ $<$ $0$ decided at $\tau_i$, but executed at times $\tau_i+mh$.  The dynamics of $Y$ is then given by
\beqs
Y_t &=& Y_0 + \sum_{\tau_i+mh\leq t} \xi_i,
\enqs
which means that discrete trading $\Delta Y_t$ $:=$ $Y_t-Y_{t^-}$ $=$ $\xi_i$ occur at times $s$ $=$ $\tau_i+mh$, $i$ $\geq$ $1$.
In absence of trading, the  cash holdings $Z$ grows deterministically at rate $r$~: $dZ_t$ $=$ $rZ_tdt$.  When a discrete trading
$\Delta Y_t$ occurs, this results in a variation of cash holdings by $\Delta Z_t$ $:=$
$Z_t-Z_{t^-}$ $=$ $- (\Delta Y_t) S_t$, from the self-financing condition. In other words, the dynamics of $Z$ is given by
\beqs
Z_t &=& Z_0 + \int_0^t r Z_u du -  \sum_{\tau_i+mh\leq t} \xi_i . S_{\tau_i+mh}.
\enqs
The wealth process is equal to $L(S_t,Y_t,Z_t)$ $=$ $Z_t+Y_tS_t$.
This financial example corresponds to the general model \reff{eqX} with $X$ $=$ $(S,Y,Z)$, $b$ $=$ $(\beta \; 0 \; r)'$,
$\sigma$ $=$ $(\gamma \; 0 \; 0)$, and
\beqs
\Gamma(s,y,z,e) &=& \left( \begin{array}{c}
                                      s \\  e \\  z - e s
                                      \end{array} \right).
\enqs
Fix now some contingent claim characterized by its payoff at time $T$~: $H(S_T)$ for some measurable function $H$. The two following hedging and valuation criteria are very popular in finance, and may be embedded in our general framework~:

\vspace{1mm}

\noindent $\bullet$ {\it Shortfall risk hedging}. The investor is looking for a trading strategy that minimizes the shortfall risk  of the  $P\&L$ between her contingent claim and her terminal wealth,
\beqs
\inf_{\alpha\in\Ac} \E\Big[ \Big( H(S_T) -  L(S_T,Y_T,Z_T) \Big)_+ \Big].
\enqs

\noindent $\bullet$ {\it Utility indifference price}. Given an utility function $U$ for the investor, an initial capital $z$ in cash, zero in stock,
and $\kappa$ $\geq$ $0$ units of contingent claims, define the expected utility under optimal trading
\beqs
V_0(z,\kappa) &=& \sup_{\alpha\in\Ac} \E\Big[ U\big(L(S_T,Y_T,Z_T) - \kappa H(S_T)\big)\Big].
\enqs
The utility  indifference ask price $\pi_a(\kappa,z)$ is the  price at which the investor is indifferent (in the sense that her expected
utility is unchanged under optimal trading) between paying nothing and not having the claim, and receiving
$\pi_a(\kappa,z)$ now to deliver $\kappa$ units of claim at time $T$. It is then defined as the solution to
\beqs
V_0(z + \pi_a(\kappa,z),\kappa) &=&  V_0(z,0).
\enqs

\subsection{Value functions}

In order to  provide an analytic characterization of the control problem  \reff{defcon},  we need  as usual  to extend the
definition of this  control problem   to general initial conditions.  However, in contrast with classical control problems without execution delay, the
diffusion process solution to \reff{eqX} is not Markovian.  Actually,  given an impulse control, we see that the state of
the system is  not only defined by its current state value at time $t$  but also by  the pending orders, that is the orders
not yet executed, i.e.  decided between time $t-mh$ and $t$. Notice that the number of pending orders is less or equal to $m$.
Let us then introduce the following definitions and notations.  For any $t$ $\in$ $[0,T]$, $k$ $=$ $0,\ldots,m$,  we denote  by
\beqs
P_t(k) &=& \Big\{  p = (t_i,e_i)_{1\leq i\leq k} \in ([0,T]\times E)^k~:   t_i - t_{i-1} \geq h, \;\; i=2,\ldots,k, \\
& &  \hspace{35mm} \;\;\;  t-mh < t_i  \leq  t,  \; i=1,\ldots,k  \Big\},
\enqs
the set of $k$ pending orders not yet executed before time $t$,
with the convention that $P_t(0)$ $=$ $\emptyset$.  For any $p$ $=$ $(t_i,e_i)_{1\leq i\leq k}$
$\in$ $P_t(k)$, $t$ $\in$ $[0,T]$, $k$ $=$ $0,\ldots,m$, we denote
\beqs
\Ac_{t,p} &=& \Big\{ \alpha = (\tau_i,\xi_i)_{i\geq 1} \in \Ac~:  (\tau_i,\xi_i) = (t_i,e_i), \; i=1,\ldots,k \;  \mbox{ and } \; \tau_{k+1} \geq t \Big\},
\enqs
the set of admissible impulse controls  with  pending orders $p$ before  time $t$.

 For any $(t,x)$ $\in$ $[0,T]\times\R^d$,  $p$ $\in$ $P_t(k)$, $k$ $=$ $0,\ldots,m$, and $\alpha$ $\in$ $\Ac_{t,p}$, we denote by
$X^{t,x,p,\alpha}$  the solution to \reff{eqX} for $t$ $\leq$ $s$ $\leq$ $T$, with initial data $X_{t}$ $=$ $x$, and pending orders $p$, i.e.
\beqs
 \hspace{-7mm} X_s  &=& x + \int_t^s  b(X_u) du + \int_t^s \sigma(X_u) dW_u +
  \sum_{t<\tau_i+mh \leq s} \big( \Gamma(X_{(\tau_i+mh)^-},\xi_i) -  X_{(\tau_i+mh)^-}\big).
\enqs
Using standard arguments based on Burkholder-Davis-Gundy's inequality, Gronwall's lemma and \reff{hypgam}, we easily check that
\beq \label{estimX2}
\E\big[ \sup_{t\leq s \leq T} |X_s^{t,x,p,\alpha}|^2\big] & \leq & C(1 + |x|^2),
\enq
for some positive constant $C$ depending only on $b$, $\sigma$, $\Gamma$ and $T$.  We then consider the following performance criterion~:
\beqs
J_k(t,x,p,\alpha) &=& \E\Big[ \int_t^T f(X_s^{t,x,p,\alpha}) ds + g(X_T^{t,x,p,\alpha}) +
\sum_{ t < \tau_i + mh \leq T} c(X_{(\tau_i+mh)^-}^{t,x,p,\alpha},\xi_i) \Big],
\enqs
for $(t,x)$ $\in$ $[0,T]\times\R^d$, $p$ $\in$ $P_t(k)$, $k$ $=$ $0,\ldots,m$, $\alpha$ $=$ $(\tau_i,\xi_i)_i$  $\in$ $\Ac_{t,p}$,
and the corresponding value functions~:
\beqs
v_k(t,x,p) &=& \sup_{\alpha \in \Ac_{t,p}}  J_k(t,x,p,\alpha), \;\;\;  \; k=0,\ldots,m, \;  (t,x,p) \in  \Dc_k,
\enqs
where $\Dc_k$ is the definition domain of $v_k$~:
\beqs
\Dc_k &=& \big\{ (t,x,p)~: (t,x) \in [0,T]\times\R^d, \; p \in P_t(k) \big\}.
\enqs
For $k$ $=$ $0$, $P_t(0)$ $=$ $\emptyset$,  and we write by convention  $v_0(t,x)$ $=$ $v_0(t,x,\emptyset)$, $\Dc_0$ $=$ $[0,T]\times\R^d$
so that the original control problem in \reff{defcon} is given by $V_0$ $=$ $v_0(0,X_0)$.
Notice  from \reff{hypfgc} and \reff{estimX2} that the functions $v_k$ satisfy the linear growth condition on $\Dc_k$~:
\beq \label{growthvk}
\sup_{(t,x,p) \in \Dc_k} \frac{|v_k(t,x,p)|}{1+|x|} & < &  \infty, \;\;\; k=0,\ldots,m.
\enq

\section{Dynamic programming}

\setcounter{equation}{0} \setcounter{Assumption}{0}
\setcounter{Theorem}{0} \setcounter{Proposition}{0}
\setcounter{Corollary}{0} \setcounter{Lemma}{0}
\setcounter{Definition}{0} \setcounter{Remark}{0}

In this section, we state the dynamic programming relation on the value functions of our control problem with delay execution.
For any $t$ $\in$ $[0,T]$, $\alpha$ $=$ $(\tau_i,\xi_i)_{i\geq 1}$ $\in$ $\Ac$, we  denote~:
 \beq
 \iota(t,\alpha) &=& \inf\{ i \geq 1~: \tau_i > t - mh \} - 1 \; \in \; \N \cup \{\infty\}, \label{defl} \\
 k(t,\alpha) & =&  {\rm card}\big\{ i \geq 1~:   t - mh < \tau_i \leq  t \big\} \; \in \;  \{0,\ldots,m\},  \label{defk} \\
 p(t,\alpha) & = &  (\tau_{i+\iota(t,\alpha)},\xi_{i+\iota(t,\alpha)})_{1\leq i\leq k(t,\alpha)}     \; \in \; P_t(k(t,\alpha)).  \label{defp}
 \enq

\begin{Theorem}
The value functions satisfy the dynamic programming principle~: for all $k$ $=$ $0,\ldots,m$, $(t,x,p)$ $\in$ $\Dc_k$,
\beq
v_k(t,x,p) &=& \sup_{\alpha\in\Ac_{t,p}} \E\Big[ \int_t^\theta f(X_s^{t,x,p,\alpha}) ds +  \sum_{\tau_i+mh \leq \theta} c(X_{(\tau_i+mh)^-}^{t,x,p,\alpha},\xi_i)  \nonumber  \\
& & \hspace{25mm} + \;  v_{k(\theta,\alpha)}(\theta,X_\theta^{t,x,p,\alpha},p(\theta,\alpha)) \Big], \label{reldynpro}
\enq
where $\theta$ is any stopping time valued in $[t,T]$, possibly depending on $\alpha$ in \reff{reldynpro}.  This means

\noindent (i) for all $\alpha$ $\in$ $\Ac_{t,p}$, for all $\theta$ stopping time valued in $[t,T]$,
\beq
v_k(t,x,p) & \geq &   \E\Big[ \int_t^\theta f(X_s^{t,x,p,\alpha}) ds +  \sum_{t<\tau_i+mh \leq \theta} c(X_{(\tau_i+mh)^-}^{t,x,p,\alpha},\xi_i) \nonumber  \\
& &   \hspace{25mm}  + \;    v_{k(\theta,\alpha)}(\theta,X_\theta^{t,x,p,\alpha},p(\theta,\alpha)) \Big]. \label{reldynpro1}
\enq
\noindent (ii) for all $\eps$ $>$ $0$, there exists $\alpha$ $\in$ $\Ac_{t,p}$ such that for all $\theta$ stopping time valued in $[t,T]$,
\beq
v_k(t,x,p) - \eps & \leq &   \E\Big[ \int_t^\theta f(X_s^{t,x,p,\alpha}) ds +  \sum_{t<\tau_i+mh \leq \theta} c(X_{(\tau_i+mh)^-}^{t,x,p,\alpha},\xi_i) \nonumber  \\
& &   \hspace{25mm}  + \;    v_{k(\theta,\alpha)}(\theta,X_\theta^{t,x,p,\alpha},p(\theta,\alpha)) \Big]. \label{reldynpro2}
\enq

\end{Theorem}

We  now give an  explicit consequence of the above dynamic programming that will be useful in the derivation of the corresponding
analytic characterization. We introduce some additional notations. For all $t$ $\in$ $[0,T]$,  we denote by $\Ic_t$ the set of pairs $(\tau,\xi)$ where  $\tau$ is a stopping time, $t$ $\leq$ $\tau$  a.s.,
and $\xi$ is a $\Fc_\tau$-measurable random variable valued in $E$. For any $p$ $=$ $(t_i,e_i)_{1\leq i\leq k}$ $\in$ $P_t(k)$, we denote $p_-$ $=$ $(t_i,e_i)_{2\leq i\leq k}$ with the convention that $p_-$ $=$ $\emptyset$ when $k$ $=$ $1$.

When no impulse control is applied to the system, we denote by $X_s^{t,x,0}$  the solution to \reff{eqX0} with initial data $X_t$ $=$ $x$, and by $\Lc$ the associated
infinitesimal generator~:
\beqs
\Lc\varphi &=& b(x).D_x\varphi + \frac{1}{2}{\rm tr}(\sigma\sigma'(x)D_x^2\varphi).
\enqs

For $k$ $\in$ $\{1,\ldots,m\}$, we partition the set $P_t(k)$ into $P_t(k)$ $=$ $P_t^1(k)$ $\cup$ $P_t^2(k)$ where
\beqs
P_t^1(k) &=& \Big\{  p = (t_i,e_i)_{1\leq i\leq k} \in P_t(k)~:   t_k > t -h   \Big\} \\
P_t^2(k) &=& \Big\{  p = (t_i,e_i)_{1\leq i\leq k}  \in P_t(k)~:  t_k \leq t-h  \Big\}.
\enqs
We easily see from the lag constraint on the pending orders that  $P_t^2(k)$ $=$ $\emptyset$ iff $k$ $=$ $m$, and so 
$P_t(m)$ $=$ $P_t^1(m)$.

%\begin{Remark}
%{\rm
%For any $k$ $=$ $1,\ldots,m$, if $p$ $=$ $(t_i,e_i)_{1\leq i\leq k}$ $\in$ $P_t(k)$, we notice from the
%lag constraint on the pending orders that $t_j-t_i$ $\geq$ $(j-i)h$, $1\leq i <j\leq k$, and $t_i$ $>$ $t-(m-i+1)h$, $i$ $=$ $1,\ldots,k$.  For $k$ $=$ $m$, %we have  $t_i$ $\in$ $(t-(m-i+1)h,t-(m-i)h]$, for all $i$ $=$ $1,\ldots,m$. In particular, $P_t^2(k)$ $=$ $\emptyset$ iff $k$ $=$ $m$, and so
%$P_t(m)$ $=$ $P_t^1(m)$.
%}
%\end{Remark}

\begin{Corollary} \label{coroldyn} Let $(t,x)$ $\in$ $[0,T)\times\R^d$.

\noindent (1) For  $k$ $\in$ $\{1,\ldots,m\}$,  and $p$  $=$ $(t_i,e_i)_{1\leq i\leq k}$ $\in$ $P_t^1(k)$ such that
$t_1+mh$ $\leq$ $T$, we have  for  any stopping time $\theta$ valued in $[t,(t_k+h)\wedge(t_1+mh))$~:
\beq \label{dynprovkpc1}
v_k(t,x,p) &=& \E\Big[ \int_t^{\theta} f(X_s^{t,x,0}) ds + v_k(\theta,X_{\theta}^{t,x,0},p) \Big].
\enq
(2)  For  $k$  $\in$ $\{0,\ldots,m-1\}$,  and $p$ $=$ $(t_i,e_i)_{1\leq i\leq k}$  $\in$ $P_t^2(k)$ such that  $t_1+mh$ $\leq$ $T$,
with the convention that $P_t^2(k)$ $=$ $\emptyset$ and $t_1+mh$ $=$ $T$ when $k$ $=$ $0$, we have
for  any stopping time $\theta$ valued in $[t,(t_1+mh)\wedge (t+h)) $~:
\beq
v_k(t,x,p) &=& \sup_{(\tau,\xi)\in\Ic_t} \E\Big[ \int_t^{\theta} f(X_s^{t,x,0}) ds +
v_k(\theta,X_\theta^{t,x,0},p)  1_{\theta < \tau}  \label{dynprovkpc2} \\
& &   \hspace{25mm} + \; v_{k+1}(\theta,X_{\theta}^{t,x,0}, p \cup (\tau,\xi))   1_{\tau \leq \theta} \Big],  \nonumber
\enq
\end{Corollary}

\vspace{2mm}

\noindent {\bf Interpretation and remarks}

\vspace{1mm}

\noindent  (1) $P_t^1(k)$ represents the set of $k$ pending orders where the last order is within the period $(t-h,t]$ of nonintervention before $t$.  Hence, from  time $t$ and until time $(t_k+h)\wedge(t_1+mh)$,  we cannot intervene on the diffusion system and no pending order will be executed  during this time period. This is mathematically formalized by relation \reff{dynprovkpc1}.

\noindent (2)  $P_t^2(k)$ represents the set of $k$ pending orders where the last order is  out of the period of nonintervention before $t$.  Hence, at time $t$,  one has  two possible decisions~: either one lets continue the system or one immediately intervene.
In this latter case,  this  order adds to the previous ones.  The mathematical formalization of these two choices is translated into relation \reff{dynprovkpc2}.

\vspace{1mm}

In the next sections, we show how one can exploit these dynamic programming relations in order to characterize analytically the value functions by
means of  partial differential equations, and then to provide an algorithm for computing the value functions.

\section{PDE system viscosity characterization}

\setcounter{equation}{0} \setcounter{Assumption}{0}
\setcounter{Theorem}{0} \setcounter{Proposition}{0}
\setcounter{Corollary}{0} \setcounter{Lemma}{0}
\setcounter{Definition}{0} \setcounter{Remark}{0}

For  $k$ $=$ $1,\ldots,m$,  let us introduce the subspace  $\Theta_k$ of $[0,T]^k$~:
\beqs
\Theta_k &=& \Big\{ t^{(k)} = (t_i)_{1\leq i\leq k} \in [0,T]^k~: t_k - t_1 < mh, \; t_i-t_{i-1} \geq h, \; i=2,\ldots,k \Big\}.
\enqs
We shall write,  by misuse of notation,  $p$ $=$ $(t_i,e_i)_{1\leq i\leq k}$ $=$ $(t^{(k)},e^{(k)})$,
for any  $t^{(k)} = (t_i)_{1\leq i\leq k} \in \Theta_k$, $e^{(k)} = (e_i)_{1\leq i\leq k} \in E^k$.  By convention, we set
$\Theta_k$ $=$ $E^k$ $=$ $\emptyset$ for $k$ $=$ $0$.
Notice that for all  $t$ $\in$ $[0,T]$, and $p$ $=$ $(t^{(k)},e^{(k)})$ $\in$ $\Theta_k\times E^k$, $k$ $=$ $0,\ldots,m$, we have
\beqs
p \in P_t(k) & \Longleftrightarrow & t \in   \T_{p}(k),
\enqs
where $\T_p(k)$ is the time domain in $[0,T]$ defined by~:
\beqs
\T_p(k) &=&  [t_k,t_1+mh) \cap [t_k,T].
\enqs
By convention, we set $\T_p(k)$ $=$ $[0,T]$ for $k$ $=$ $0$.
We can then rewrite the domain $\Dc_k$ of the value function $v_k$ in terms of union of time-space domains~:
\beqs
\Dc_k  &=&   \Big\{ (t,x,p )~:  (t,x) \in \T_p(k) \times \R^d, \; p   \in   \Theta_k \times E^k \Big\}.
\enqs
Therefore, the determination of the value function $v_k$, $k$ $=$ $0,\ldots,m$,  is  equivalent to the determination of the
function  $v_k(.,.,p)$ on $\T_p(k)\times\R^d$ for all $p$  $\in$ $\Theta_k\times E_k$.
The main goal of this paper is to provide an analytic characterization of these functions by means of the dynamic programming principle stated in the previous section.

For $k$ $=$ $1,\ldots,m$,   we denote
\beqs
\Theta_k^m &=& \Big\{ t^{(k)} = (t_i)_{1\leq i\leq k} \in \Theta_k~:  t_1 + mh  \leq T  \Big\},  \;\;\;  \Theta_k(m) \; = \;  \Theta_k\setminus\Theta_k^m,
\enqs
and we define the ``$m$-interior" of $\Dc_k$ by~:
\beqs
\Dc_k^m &=& \Big\{ (t,x,p) \in \Dc_k~: p \in \Theta_k^m \times E^k \Big\}.
\enqs
For $k$ $=$ $0$, we set  $\Dc_0^m$ $=$ $[0,T)\times\R^d$.
For $p$ $=$ $(t_i,e_i)_{1\leq i\leq k}$  $\in$ $\Theta_k^m\times E^k$,
we partition the  time domain $\T_p(k)$ into  $\T_p(k)$ $=$ $\T_p^1(k)\cup \T_p^2(k)$ where
\beqs
\T_p^1(k) &=& \Big\{ t \in \T_p(k)~:  t < t_k + h\Big\} \; = \; [t_k,  (t_k+h)\wedge (t_1+mh))    \\
\T_p^2(k) &=& \Big\{ t \in \T_p(k)~:   t \geq t_k + h     \Big\}  \; = \; [t_k+h,t_1+mh),
\enqs
with the convention that $[s,t)$ $=$ $\emptyset$ if $s$ $\geq$ $t$.  We then partition $\Dc_k^m$ into
$\Dc_k^m$ $=$ $\Dc_k^{1,m}$ $\cup$ $\Dc_k^{2,m}$ where
\beqs
\Dc_k^{1,m} &=& \big\{ (t,x,p) \in \Dc_k^m~: t \in \T_p^1(k) \Big\} \\
\Dc_k^{2,m} &=& \big\{ (t,x,p) \in \Dc_k^m~: t \in \T_p^2(k) \Big\}.
\enqs
Notice that  for $k$ $=$ $1,\ldots,m$, and any $p$ $\in$ $\Theta_k^m\times E^k$,  $\T_p^1(k)$ is never empty. In particular,
$\Dc_k^{1,m}$ $\neq$ $\emptyset$.  For $k$ $=$ $m$, and any $p$ $=$ $(t_i,e_i)_{1\leq i\leq m}$
$\in$ $\Theta_m\times E^m$,  we have $t_m+h$ $\geq$ $t_1+mh$, and so $\T_p^2(m)$ $=$ $\emptyset$.  Hence, $\Dc_m^{2,m}$ $=$ $\emptyset$ and
$\Dc_m^m$ $=$ $\Dc_m^{1,m}$.

\vspace{2mm}

The PDE system to our control problem is formally derived by sending $\theta$ to
$t$ $<$ $t_1+mh$  into  dynamic programming relations \reff{dynprovkpc1}-\reff{dynprovkpc2}. This  provides  equations for the value functions $v_k$ on $\Dc_k^m$, which take the following nonstandard form, and are  divided into~:
\beq \label{viscoQ1}
 - \Dt{v_k}(t,x,p)  - \Lc v_k(t,x,p)  - f(x)   &=& 0 \;\;\;   \mbox{ on  } \;\;\;  \Dc_k^{1,m},   \;\; \; k=1,\ldots,m,
\enq
\beq
 & & \hspace{-17mm} \min\Big\{ - \Dt{v_k}(t,x,p)  - \Lc v_k(t,x,p) - f(x)  \; ,   \nonumber  \\
& &    v_k(t,x,p)  -   \sup_{e \in E}  v_{k+1}(t,x,p\cup (t,e)) \Big\} \; = \;  0 \;\;\;  \mbox{ on }  \;\;\;  \Dc_k^{2,m}, \;\; k=0,\ldots,m-1,   \label{viscoQ2}
\enq
with the convention that $\Dc_0^{2,m}$ $=$ $\Dc_0^m$ $=$ $[0,T)\times\R^d$.

 \vspace{2mm}

As usual, the value functions  need  not be  smooth, and even not known to be continuous a priori, and we shall work with
the notion of (discontinuous) viscosity solutions (see \cite{craishlio92} or \cite{fleson93} for  classical references on the subject), which we adapt in our context as follows. For a locally bounded function $w_k$ on $\Dc_k^m$, we denote $\underline{w_k}$ (resp. $\overline w_k$) its lower semicontinuous
(resp. upper-semicontinuous) envelope, i.e.
\beqs
\underline{w_k}(t,x,p) & = &    \liminf_{\tiny{(t',x',p')  \rightarrow (t,x,p)}} w_k(t',x',p'), \\
 \overline{w_k}(t,x,p)  &= & \limsup_{\tiny{(t',x',p')  \rightarrow (t,x,p)}} w_k(t',x',p'), \;\;\; (t,x,p) \in \Dc_k^m, \; k=0,\ldots,m.
\enqs

\begin{Definition}
We say that  a family of locally bounded  functions $w_k$ on $\Dc_k^m$, $k$ $=$ $0,\ldots,m$,  is a viscosity supersolution  (resp. subsolution)  of  \reff{viscoQ1}-\reff{viscoQ2} on $\Dc_k^m$, $k$ $=$ $0,\ldots,m$,  if~:

\noindent (i)  for all $k$ $=$ $1,\ldots,m$,  $(t_0,x_0,p_0)$ $\in$ $\Dc_k^{1,m}$, and $\varphi$ $\in$
$C^{2}(\Dc_k^{1,m})$, which realizes a local minimum of $\underline{w_k}-\varphi$  (resp. maximum of
$\overline{w_k}-\varphi$), we have
\beqs
 - \Dt{\varphi}(t_0,x_0,p_0)  - \Lc \varphi(t_0,x_0)  - f(x_0)   & \geq  & 0  \;\;\; (resp. \; \leq \; 0).
\enqs
(ii) for all $k$ $=$ $0,\ldots,m-1$,   $(t_0,x_0,p_0)$ $\in$ $\Dc_k^{2,m}$, and $\varphi$ $\in$
$C^{2}(\Dc_k^{2,m})$, which realizes a local minimum of $\underline{w_k}-\varphi$  (resp. maximum of
$\overline{w_k}-\varphi$), we have
\beqs
 \min\big\{  - \Dt{\varphi}(t_0,x_0,p_0)  - \Lc \varphi(t_0,x_0,p_0)  - f(x_0) \; , \; \hspace{7mm}  & &   \\
\;\;\;\;\;\;\;   \underline{w_k}(t_0,x_0,p_0) -  \sup_{e\in E} \underline{w_{k+1}}(t_0,x_0,p_0 \cup (t_0,e)) \big\} &  \geq &   0
\enqs
(resp.
\beqs
 \min\big\{  - \Dt{\varphi}(t_0,x_0,p_0)  - \Lc \varphi(t_0,x_0,p_0)  - f(x_0) \; , \; \hspace{7mm}  & &   \\
\;\;\;\;\;\;\;   \overline{w_k}(t_0,x_0,p_0) -  \sup_{e\in E} \overline{w_{k+1}}(t_0,x_0,p_0 \cup (t_0,e)) \big\} &  \leq &   0 )
\enqs
We say that a family of locally bounded  functions $w_k$ on $\Dc_k^m$, $k$ $=$ $0,\ldots,m$,
is a viscosity solution of  \reff{viscoQ1}-\reff{viscoQ2} if it is a
viscosity supersolution and subsolution of  \reff{viscoQ1}-\reff{viscoQ2}.
\end{Definition}

\vspace{2mm}

We then state the viscosity  property of the value functions to our control problem.

\begin{Proposition} \label{provisco} (Viscosity property)

\noindent  The family of value functions $v_k$, $k$ $=$ $0,\ldots,m$,  is a viscosity solution to \reff{viscoQ1}-\reff{viscoQ2}. Moreover,
for all $k$ $=$ $0,\ldots,m-1$,  $(t,x,p)\in \Dc_{k}^{m}$, $p$ $=$ $(t_i,e_i)_{1\leq i\leq k}$ with  $t=t_{k}+h$,  we have~:
\beq \label{condSupplDisc}
\underline{v_{k}}(t,x,p)\geq \sup_{e\in E} \underline{v_{k+1}}(t,x,p\cup(t,e)).
\enq
\end{Proposition}

In order to have a complete characterization of the value functions, and so of our control  problem, we need to determine the suitable boundary conditions. These concern for $k$ $=$ $1,\ldots,m$ the time-boundary of $\Dc_k^m$, i.e. the points $(t_1+mh,x,p)$ for $x$ $\in$ $\R^d$, $p$ $=$ $(t_i,e_i)_{1\leq i\leq k}$ $\in$ $\Theta_k^m\times E^k$, and also the complement set of $\Dc_k^m$ in $\Dc_k$.
For a locally bounded function $w_k$ on $\Dc_k^m$, $k$ $=$ $1,\ldots,m$, we denote
\beqs
\overline{w_k}(t_1+mh,x,p)  &= &  \limsup_{\tiny{\begin{array}{c}(t,x',p')\rightarrow (t_1+mh,x,p) \\ (t,x',p')\in\Dc_k^m \end{array}}} w_k(t,x',p'), \\
\underline{w_k}(t_1+mh,x,p) & = &  \liminf_{\tiny{\begin{array}{c}(t,x',p')\rightarrow (t_1+mh,x,p) \\ (t,x',p')\in\Dc_k^m \end{array}}} v_k(t,x',p'), \;\; x \in \R^d, \; p =(t_i,e_i)_{1\leq i\leq k}  \in \Theta_k^m,
\enqs
and if  these two limits are equal, we set
\beqs
w_k((t_1+mh)^-,x,p) & = &  \overline{w_k}(t_1+mh,x,p)  \; = \;  \underline{w_k}(t_1+mh,x,p).
\enqs
We also denote for a locally bounded function $w_0$ on $[0,T)\times\R^d$~:
\beqs
\overline{w_0}(T,x) \; = \;  \limsup_{t\nearrow T, x'\rightarrow x} w_0(t,x'), & &
\underline{w_0}(T,x) \; = \;  \liminf_{t\nearrow T, x'\rightarrow x} w_0(t,x'), \;\;\; x \in \R^d,
\enqs
and if these two limits are equal, we set $w_0(T^-,x)$ $=$ $\overline{w_0}(T,x)$ $=$ $\underline{w_0}(T,x)$.
The complement set of $\Dc_k^m$ in $\Dc_k$ is
\beqs
\Dc_k(m) \; := \; \Dc_k \setminus\Dc_k^m &=& \big\{ (t,x,p) \in \Dc_k~:  p \in \Theta_k(m)\times E^k\Big\}.
\enqs

\begin{Proposition} \label{propdata} (Boundary data)

\noindent (i) For $k$ $=$ $1,\ldots,m$,  $p$ $=$ $(t_i,e_i)_{1\leq i\leq k}$ $\in$ $\Theta_k^m\times E^k$, $x$ $\in$ $\R^d$, $v_k((t_1+mh)^-,x,p)$ exists and~:
\beq
v_k((t_1+mh)^-,x,p) &=& c(x,e_1) + v_{k-1}(t_1+mh,\Gamma(x,e_1),p_-).   \label{vkvk-1}
\enq
(ii) For $k$  $=$ $1,\ldots,m$,  we have~:
\beq \label{vkT}
v_k(t,x,p) &=& \E\Big[ \int_t^T f(X_s^{t,x,0}) ds + g(X_T^{t,x,0}) \Big], \;\;\; (t,x,p) \in \Dc_k(m).
\enq
\end{Proposition}

\vspace{2mm}

We can now state the unique PDE characterization result for our control delay problem.

\begin{Theorem} \label{thmmain}
The family of value functions $v_k$, $k$ $=$ $0,\ldots,m$, is the unique viscosity solution to \reff{viscoQ1}-\reff{viscoQ2}, which satisfy \reff{condSupplDisc}, the boundary data \reff{vkvk-1}-\reff{vkT}, and the linear growth condition \reff{growthvk}.
Moreover, $v_k$ is continuous on $\Dc_k^m$ and on $\Dc_k(m)$, $k$ $=$ $0,\ldots,m$.
\end{Theorem}

\begin{Remark} {\rm  {\bf (Case $m$ $=$ $1$)}

\noindent  In  the particular case where the execution delay is equal to the intervention lag, i.e. $m$ $=$ $1$, we have  two value functions $v_0$ and
$v_1$, and  the system \reff{viscoQ1}-\reff{viscoQ2} may be  significantly  simplified.  Actually,
from the linear PDE \reff{viscoQ1} and the boundary data \reff{vkvk-1} for $k$ $=$ $m$ $=$ $1$, we have the Feynman-Kac representation~:
\beq \label{v0v1}
v_1(t,x,(t_1,e_1)) &=& \E\Big[ \int_t^{t_1+h} f(X_s^{t,x,0}) ds + c( X_{t_1+h}^{t,x,0},e)
+ v_{0}(t_1+h,\Gamma(X_{t_1+h}^{t,x,0},e)) \Big],
\enq
for all  $(t_1,e_1)$ $\in$ $[0,T-h]\times E$,  $(t,x)$ $\in$ $[t_1,t_1+h)\times\R^d$.
By plugging \reff{v0v1} for $t$ $=$ $t_1$ into \reff{viscoQ2} for $k$ $=$ $0$, we obtain the variational inequality satisfied by $v_0$~:
\beq
 0 \; = \; \min\Big\{ - \Dt{v_0} - \Lc v_0 - f \; ,  \hspace{35mm} & &  \label{inegv0} \\
v_0 -   \sup_{e \in E}   \E\Big[ \int_t^{t+h} f(X_s^{t,x,0}) ds + c( X_{t+h}^{t,x,0},e)
+ v_{0}(t+h,\Gamma(X_{t+h}^{t,x,0},e)) \Big] \Big\} & \mbox{on} &  [0,T-h] \times\R^d, \nonumber
\enq
together with the terminal condition for $k$ $=$ $0$ (see \reff{v0f0})~:
\beq \label{limv0}
v_0(t,x) &=& \E\Big[ \int_t^T f(X_s^{t,x,0}) ds + g(X_T^{t,x,0}) \Big],  \;\;\; (t,x) \in (T-h,T] \times\R^d.
\enq
Therefore, in the case $m$ $=$ $1$,  and as observed in \cite{rob76},  the original problem  is reduced to a no-delay impulse control problem
\reff{inegv0} for $v_0$, and  $v_1$ is explicitly related to $v_0$ by \reff{v0v1}.  Equations \reff{inegv0}-\reff{limv0} can be solved by  iterated optimal stopping problems, see the details in the next section in the more general case $m$ $\geq$ $1$.
 }
 \end{Remark}

\begin{Remark}
{\rm In the general case $m$ $\geq$ $1$, we point out the peculiarities of the PDE characterization for our control delay problem.

\noindent {\bf 1.} The dynamic programming coupled system \reff{viscoQ1}-\reff{viscoQ2} has a nonstandard form. For fixed $k$,
there is a discontinuity on the differential operator of  the equation satisfied by $v_k$ on $\Dc_k^m$.  Indeed, the PDE is divided into a linear equation on the subdomain $\Dc_k^{1,m}$, and a variational inequality  with obstacle involving the value function $v_{k+1}$ on the
subdomain $\Dc_k^{2,m}$.  Moreover,  the time domain  $\T_p(k)$ of $\Dc_k^m$ for $v_k(.,x,p)$ depends on the argument $p$ $\in$ $\Theta_k^m$.
With respect to usual comparison principle of nonlinear PDE,  we state an uniqueness result  for viscosity solutions
satisfying in addition the inequality \reff{condSupplDisc} at the discontinuity of the differential operator.

\noindent {\bf 2.}  The  boundary data also present some specificities.  For fixed $k$, the condition in \reff{vkvk-1} concerns  as usual data on the time-boundary of the domain $\Dc_k^m$ on which the value function $v_k$ satisfies a PDE.  However,  it involves
data on the value function $v_{k-1}$, which is a priori not known.
The condition in \reff{vkT} for $v_k$  concerns the complement set of $\Dc_k^m$, and is explicitly known.
Notice also  that we do not need to specify in Theorem \ref{thmmain} the boundary data for $v_0$. Actually, this will be derived in \reff{v0f0} as a direct consequence of \reff{vkT} for $k$ $=$ $1$ and the PDE equation \reff{viscoQ2} for $k$ $=$ $0$.

\noindent {\bf 3.} The continuity property of the value functions $v_k$ on $\Dc_k^m$ is not at all obvious a priori from the very
definitions of $v_k$, and is proved  actually as consequences of  comparison principles and boundary data for
the system \reff{viscoQ1}-\reff{viscoQ2}, see Proposition \ref{procont}. The continuity of $v_k$ on $\Dc_k(m)$ is obvious from the boundary data \reff{vkT}. We mention, however, that the value functions $v_k$, $k$ $\geq$ $1$, are not continuous in general on their whole domain $\Dc_k$~: there is a discontinuity at points $(T,x,p)$ with $p$ $=$ $(t_i,e_i)_{1\leq i\leq k}$ $\in$ $\Theta_k$ s.t. $t_1+mh$ $=$ $T$. Indeed, from the very definition of the value functions, we have for such points $v_k(T,x,p)$ $=$ $g(x)$ (and also  $v_0(T,x)$ $=$ $g(x)$), while from \reff{vkvk-1}, we have
$v_k(T^-,x,p)$ $=$ $c(x,e_1)$ $+$ $g(x)$. Hence, $v_k(T,x,p)$ $\neq$ $v_k(T^-,x,p)$ once $c(x,e_1)$ $\neq$ $0$.
 }
 \end{Remark}

\vspace{3mm}

The  PDE characterization in Theorem \ref{thmmain} means that the value functions are in theory completely determined by the resolution of the PDE system \reff{viscoQ1}-\reff{viscoQ2} together with the boundary data \reff{vkvk-1}-\reff{vkT}. We show in the next section how to solve this system and compute in practice these value functions and the associated optimal impulse controls.

\section{An algorithm to compute the value functions and the optimal control} \label{secalgo}

\setcounter{equation}{0} \setcounter{Assumption}{0}
\setcounter{Theorem}{0} \setcounter{Proposition}{0}
\setcounter{Corollary}{0} \setcounter{Lemma}{0}
\setcounter{Definition}{0} \setcounter{Remark}{0}

\subsection{Computation of the value functions}

We first make the following observation. Let us denote by $F_0$ the function defined on $[0,T]\times\R^d$ by
\beqs
F_0(t,x) &=& \sup_{e\in E} v_1(t,x,(t,e)).
\enqs
From \reff{vkT} for $k$ $=$ $1$, we deduce that for all  $e$ $\in$ $E$,
\beqs
F_0(t,x) \; = \; v_1(t,x,(t,e)) &=& \E\Big[ \int_t^T f(X_s^{t,x,0}) ds + g(X_T^{t,x,0}) \Big], \;\; (t,x) \in (T-mh,T]\times\R^d.
\enqs
This function $F_0$ clearly satisfies the linear PDE~: $-\Dt{F_0} - \Lc F_0 - f$ $=$ $0$.
Hence,  with \reff{viscoQ2} for $k$ $=$ $0$, this shows that
\beq \label{v0f0}
v_0(t,x) &=& F_0(t,x),  \;\;\; \; \; (t,x) \in (T-mh,T]\times\R^d,
\enq
and in particular, $v_0(T^-,x)$ $=$ $F_0(T,x)$ $=$ $g(x)$.  Together with  the PDE \reff{viscoQ2} for $k$ $=$ $0$, and a standard uniqueness result for  the corresponding free-boundary problem,  this proves that $v_0$ may also be represented as the solution to the optimal stopping problem~:
\beq \label{v0opt}
v_0(t,x) &=& \sup_{\tau \in \Tc_{t,T}} \E[ F_0(\tau,X_\tau^{t,x,0}) ], \;\;\; (t,x) \in [0,T]\times\R^d,
\enq
where $\Tc_{t,T}$ denotes the set of stopping times $\tau$ valued in $[t,T]$.  Hence, the value function $v_0$ is completely determined once we can compute $v_1$.

We show how one can compute
$v_k(.,.,p)$ on $\T_p(k)\times\R^d$ for all $p$ $\in$ $\Theta_k\times E^k$,  $k$ $=$ $1,\ldots,m$ and $v_0$ on $[0,T]\times\R^d$.

\vspace{2mm}

\noindent For $k$ $=$ $1,\ldots,m$, and any $n$ $\geq$ $1$, we denote~:
\beqs
\Theta_k(n) &=& \Big\{  t^{(k)} = (t_i)_{1\leq i\leq k} \in \Theta_k~:  t_1 > T -nh \Big\}, \\
N &=&  \inf\{ n \geq 1~: T -nh  < 0 \},
\enqs
so that  $\Theta_k(n)$ is strictly included in $\Theta_k(n+1)$ for $n$ $=$ $1,\ldots,N-1$, and $\Theta_k(N)$ $=$ $\Theta_k$.
We also denote for $k$ $=$ $0$, and $n$ $\geq$ $1$, $\T^n(0)$ $=$ $(T-nh,T]$ $\cap$ $[0,T]$ so that $\T^n(0)$ $=$ $(T-nh,T]$ is increasing with $n$ $=$ $1,\ldots,N-1$, and $\T^N(0)$ $=$ $[0,T]$.   We assumed $T-mh$ $\geq$ $0$ to avoid trivialities so that  $N$ $>$ $m$.
We denote for $k$ $=$ $0,\ldots,m$, and $n$ $=$ $m,\ldots,N$, 
\beqs
\Dc_k(n) &=& \Big\{ (t,x,p) \in \Dc_k~:  p \in \Theta_k(n)\times E^k \Big\}, \\
\Dc_k^m(n) & = & \Dc_k(n) \cap \Dc_k^m \; = \; \Big\{ (t,x,p) \in \Dc_k~:  p \in \Theta_k^m(n)\times E^k \Big\}, \;
\Theta_k^m(n)  =  \Theta_k(n)\setminus\Theta_k(m) \\
\Dc_k^{i,m}(n) & = &  \Dc_k(n) \cap \Dc_k^{i,m}  \; = \; \Big\{ (t,x,p) \in \Dc_k^m(n)~:  t \in \T_p^i(k)  \Big\}, \;\; i=1,2,
\enqs
with the convention that $\Dc_0(n)$  $=$ $\T^n(0)\times\R^d$, so that $\Dc_k(n)$ is strictly included in $\Dc_k(n+1)$ for
$n$ $=$ $1,\ldots,N-1$, and $\Dc_k(N)$ $=$ $\Dc_k$.
We shall compute  $v_k$ on $\Dc_k(n)$,  $k$ $=$ $0,\ldots,m$, by forward induction on $n$ $=$ $m,\ldots,N$ and backward induction on $k$.

\vspace{1mm}

\noindent $\tri$ {\bf Initialization phase~: $n$ $=$ $m$.}  From \reff{vkT} and \reff{v0f0}, we know the values of $v_k$ on $\Dc_k(m)$,
$k$ $=$ $0,\ldots,m$~:
\beqs
v_k(t,x,p) &=& \E\Big[ \int_t^T f(X_s^{t,x,0}) ds + g(X_T^{t,x,0}) \Big].
\enqs

\vspace{1mm}

\noindent $\tri$ {\bf Step $n$  $\rightarrow$  $n+1$ for $n$ $\in$ $\{m,\ldots,N-1\}$}.  Suppose we know the values of $v_k$ on $\Dc_k(n)$,
$k$ $=$ $0,\ldots,m$. In order to determine $v_k$ on $\Dc_k(n+1)$, $k$ $=$ $0,\ldots,m$, it suffices to compute
$v_k(.,.,p)$ on  $\T_p(k)\times\R^d$ for all $p$ $\in$ $\Theta_k^m(n+1)\times E^k$, $k$ $=$ $1,\ldots,m$,
and $v_0$ on $\T^{n+1}(0)\times\R^d$. We shall argue  by backward induction on $k$ $=$ $m,\ldots,0$.
\begin{itemize}

\item  Let $k$ $=$ $m$, and take some arbitrary  $p$ $=$ $(t_i,e_i)_{1\leq i\leq m}$ $\in$ $\Theta_{m}^m(n+1)\times E^{m}$.
Recall that $\T_p^2(m)$  is always empty so that $\T_p(m)$ $=$ $\T_p^1(m)$ $=$ $[t_m,t_1+mh)$.
From \reff{vkvk-1} for $k$ $=$ $m$, we have $v_{m}((t_1+mh)^-,x,p)$ $=$
$c(x,e_1) + v_{m-1}(t_1+mh,\Gamma(x,e_1),p_-)$ for all $x$ $\in$ $\R^d$, which is known from step $n$ since  either
$p_-$ $\in$ $\Theta_{m-1}(n)\times E^{m-1}$ when $m$ $>$ $1$, or  $t_1+mh$ $\in$ $\T^n(0)$ when  $m-1$ $=$ $0$.
We then solve  $v_{m}(.,.,p)$ on  $\T_p^1(m)\times\R^d$ from \reff{viscoQ1} for $k$ $=$ $m$,
which gives~:
\beqs
v_{m}(t,x,p) &=& \E\Big[ \int_t^{t_1+mh} f(X_s^{t,x,0}) ds + c(X_{t_1+mh}^{t,x,0},e_1) \\
& & \hspace{13mm} + \;   v_{m-1}(t_1+mh,\Gamma(X_{t_1+mh}^{t,x,0},e_1),p_-)\Big].
\enqs
We have then computed the value of $v_{m}(.,.,p)$ on $\T_p(m)\times\R^d$.

\item From $k+1$ $\rightarrow$ $k$ for $k$ $=$ $m-1,\ldots,1$.
(This step is  empty when $m$ $=$ $1$).  Suppose we know the values of $v_{k+1}(.,.,p)$ on
$\T_p(k+1)\times\R^d$ for all $p$ $\in$  $\Theta_{k+1}^m(n+1)\times E^{k+1}$.  Take now some arbitrary $p$ $=$
$(t_i,e_i)_{1\leq i\leq k}$  $\in$  $\Theta_k^m(n+1)\times E^k$.
We shall  compute  $v_k(.,.,p)$  successively on $\T_p^2(k)\times\R^d$ (if it is not empty) and
then on  $\T_p^1(k)\times\R^d$, and we distinguish the two cases~:

{\bf (i)} $\T_p^2(k)$ $=$ $\emptyset$.  This means $t_k+h$ $\geq$ $t_1+mh$, and so $\T_p(k)$ $=$ $\T_p^1(k)$ $=$ $[t_k,t_1+mh)$.  We then
compute $v_k(.,.,p)$  on  $\T_p(k)\times\R^d$ as above for $k$ $=$ $m$~:
\beqs
v_{k}(t,x,p) &=& \E\Big[ \int_t^{t_1+mh} f(X_s^{t,x,0}) ds + c(X_{t_1+mh}^{t,x,0},e_1) \\
& & \hspace{13mm} + \;   v_{k-1}(t_1+mh,\Gamma(X_{t_1+mh}^{t,x,0},e_1),p_-)\Big],
\enqs
where the r.h.s. is known from step  $n$ since either $p_-$ $\in$ $\Theta_{k-1}(n)\times E^{k-1}$  when $k$ $>$ $1$,  or
$t_1+mh$ $\in$ $\T^n(0)$ when $k-1$ $=$ $0$.

{\bf (ii)}  $\T_p^2(k)$ $\neq$ $\emptyset$.  This means $t_k+h$ $<$ $t_1+mh$, and so  $\T_p^1(k)$ $=$
$[t_k,t_k+h)$, $\T_p^2(k)$ $=$ $[t_k+h,t_1+mh)$.  For all $(t,x)$ $\in$ $\T_p^2(k)\times\R^d$, and $e$ $\in$ $E$, we have
$p'$ $=$ $p\cup (t,e)$ $\in$ $\Theta_{k+1}^m(n+1)\times E^{k+1}$, and $(t,x)$ $\in$ $\T_{p'}(k+1)\times\R^d$.
Hence, from  the induction hypothesis at order $k+1$, we know the value of the function~:
\beqs
F_{k,p}(t,x) &=& \sup_{e\in E} v_{k+1}(t,x,p\cup (t,e)), \;\;\; (t,x) \in \T_p^2(k)\times\R^d.
\enqs
We also know from step $n$ the value of the function~:
\beqs
G_{k,p}(x) &=&  c(x,e_1) + v_{k-1}(t_1+mh,\Gamma(x,e_1),p_-), \;\;\; x \in \R^d.
\enqs
Then, from  the PDE \reff{viscoQ2} and the terminal condition \reff{vkvk-1} at $k$,  we compute $v_k(.,.,p)$ on $\T_p^2(k)\times\R^d$ as the solution to
an optimal stopping problem with obstacle $F_{k,p}$ and terminal condition $G_{k,p}$~:
\beqs
v_k(t,x,p) &=& \sup_{\tau \in \Tc_{t,t_1+mh}} \E[ F_{k,p}(\tau,X_\tau^{t,x,0}) 1_{\tau < t_1+mh}  \nonumber \\
& & \;\;\;\;\; \hspace{13mm} + \;  G_{k,p}(X_{t_1+mh}^{t,x,0}) 1_{\tau=t_1+mh} ], \;\;\; (t,x) \in \T_p^2(k)\times\R^d.
\enqs
In particular, by continuity of $v_k(.,.,p)$ on $\T_p(k)$, we know the value of $\lim_{t\nearrow t_k+h} v_k(t,x,p)$ $=$ $v_k(t_k+h,p)$.
We then compute $v_k(.,.,p)$ on  $\T_p^1(k)\times\R^d$  from \reff{viscoQ1}~:
\beqs
v_{k}(t,x,p) &=& \E\Big[ \int_t^{t_k+h} f(X_s^{t,x,0}) ds + v_{k}(t_k+h,X_{t_k+h}^{t,x,0},p)\Big].
\enqs
We have then computed the value of $v_{k}(.,.,p)$ on $\T_p(k)\times\R^d$.

\item From $k$ $=$ $1$ $\rightarrow$ $k$ $=$ $0$.  From the above item, we know the value of $v_1(.,.,p)$ on $\T_p(1)\times\R^d$ for all $p$ $\in$ $\Theta_1(n+1)\times E$. Hence, we know the value of~:
\beqs
F_0(t,x) &=& \sup_{e\in E} v_1(t,x,(t,e)), \;\;\;\; \forall \; (t,x) \in \T^{n+1}(0)\times\R^d.
\enqs
From \reff{v0opt}, we then compute $v_0$ on $\T^{n+1}(0)\times\R^d$ as an optimal stopping problem with obstacle $F_1$.

\vspace{1mm}

We have then calculated $v_k(.,.,p)$ on  $\T_p(k)\times\R^d$ for all $p$ $\in$ $\Theta_k^m(n+1)\times E^k$
and $v_0$ on $\T^{n+1}(0)\times\R^d$, and step $n+1$ is  stated. Finally, at step $n$ $=$ $N$, the computation of the value functions is completed since
$\Dc_k(N)$ $=$ $\Dc_k$, $k$ $=$ $0,\ldots,m$.
\end{itemize}

\subsection{Description of the optimal impulse control}

In view of the above dynamic programming relations, and the general theory of optimal stopping (see \cite{elk79}), we can describe the structure of the optimal impulse control  for $V_0$ $=$ $v_0(0,X_0)$ in terms of the value functions.  Let us
define the following quantities~:

\vspace{1mm}

\noindent $\tri$  {\bf Initialization~: $n$ $=$ $0$}
\begin{itemize}
\item  {\it given an initial pending order number $k$ $=$ $0$}, we define
\beqs
\tilde\tau_1^{(0)} &=& \inf\big\{ t \geq 0~: v_0(t,X_t^{\alpha^*}) \; = \;  \sup_{e \in E} v_1(t,X_t^{\alpha^*},(t,e)) \big\} \wedge T, \\
\tilde e_1^{(0)} & \in & {\rm arg}\max_{e\in E} v_1(\tilde \tau_1^{(0)},X_{\tilde\tau_1^{(0)}}^{\alpha^*},(\tilde\tau_1^{(0)},e)).
\enqs
If $\tilde\tau_1^{(0)} + mh$ $>$ $T$, we stop the induction at $n$ $=$ $0$, otherwise  continue to the next item~:
\item {\it  Pending orders number $k$ $\rightarrow$ $k+1$} (this step is empty when $m$ $=$ $1$)  from $k$ $=$ $1$~:
\beqs
\tilde\tau_{k+1}^{(0)} &=& \inf\big\{ t \geq \tilde\tau_{k}^{(0)} +h ~:  \\
& & \;\;\;\;\; v_k(t,X_t^{\alpha^*}) \; = \;   \sup_{e \in E} v_{k+1}(t,X_t^{\alpha^*}, (\tilde\tau_i^{(0)},\tilde e_i^{(0)})_{_{1\leq i\leq k}} \cup  (t,e)) \big\} \wedge T, \\
\tilde e_{k+1}^{(0)} & \in & {\rm arg}\max_{e\in E} v_{k+1}(\tilde \tau_{k+1}^{(0)},X_{\tilde\tau_{k+1}^{(0)}}^{\alpha^*},
(\tilde\tau_i^{(n)},\tilde e_i^{(0)})_{_{1\leq i\leq k}} \cup (\tilde\tau_{k+1}^{(0)},e)).
\enqs
As long as  $\tilde\tau_{k}^{(0)}$ $\leq$ $\tilde\tau_{1}^{(0)}+mh$,  increment $k$ $\rightarrow$ $k+1$~: $\tilde\tau_{k}^{(0)}$
$\rightarrow$ $\tilde\tau_{k+1}^{(0)}$, until
\beqs
k_0 &=& \sup\big\{ k~:  \tilde\tau_{k}^{(0)} \leq \tilde\tau_{1}^{(0)}+mh \} \;\; \in \; \{1,\ldots,m\},
\enqs
and  increment  the induction on $n$  by  the following step ~:
\end{itemize}

\noindent $\tri$ {\bf  $n$ $\rightarrow$ $n+1$~:}
\begin{itemize}
\item  {\it given an initial  pending orders number  $k$ $=$ $k_n-1$}, we define
\beqs
\tilde\tau_{k_n}^{(n+1)} &=& \inf\big\{ t \geq (\tilde\tau_{1}^{(n)}+mh)\vee ( \tilde\tau_{k_n}^{(n)} +h) ~:  \\
& & \;\;  v_{_{k_n-1}}(t,X_t^{\alpha^*},\tilde p_{n^-}) \; = \;
\sup_{e \in E} v_{_{k_n}}(t,X_t^{\alpha^*},\tilde p_{n^-} \cup(t,e)) \big\} \wedge T, \\
\tilde e_{k_n}^{(n+1)} & \in & {\rm arg}\max_{e\in E} v_{_{k_n}}(\tilde\tau_1^{n+1},X_{\tilde\tau_1^{n+1}}^{\alpha^*},\tilde p_{n^-}
\cup(\tilde\tau_{k_n}^{n+1},e)),
\enqs
where we set $\tilde p_{n^-}$ $=$ $(\tilde\tau_i^{(n)},\tilde e_i^{(n)})_{2\leq i\leq k_n}$.  We denote $\tilde\tau_1^{(n+1)}$ $=$
$\tilde\tau_2^{(n)}$ if $k_n$ $>$ $1$, and $\tilde\tau_1^{(n+1)}$ $=$ $\tilde\tau_{k_n}^{n+1}$ if $k_n$ $=$ $1$.
If $\tilde\tau_1^{(n+1)} + mh$ $>$ $T$, we stop the induction at $n+1$, otherwise  continue to the next item~:

\item {\it  Pending orders number $k$ $\rightarrow$ $k+1$} (this step is empty when $m$ $=$ $1$) from $k$ $=$ $k_n$~:
\beqs
\tilde\tau_{k+1}^{(n+1)} &=& \inf\big\{ t \geq \tilde\tau_{k}^{(n+1)} +h ~:  \\
& & \;  v_k(t,X_t^{\alpha^*}) \; = \;   \sup_{e \in E} v_{k+1}(t,X_t^{\alpha^*},
\tilde p_{n^-} \cup (\tilde\tau_i^{(n+1)},\tilde e_i^{(n+1)})_{_{k_n\leq i\leq k}} \cup  (t,e)) \big\} \wedge T \\
\tilde e_{k+1}^{(n+1)} & \in & {\rm arg}\max_{e\in E} v_{k+1}(\tilde \tau_{k+1}^{(n+1)},X_{\tilde\tau_{k+1}^{(n+1)}}^{\alpha^*},
\tilde p_{n^-} \cup (\tilde\tau_i^{(n+1)},\tilde e_i^{(n+1)})_{_{k_n\leq i\leq k}}   \cup (\tilde\tau_{k+1}^{(n+1)},e))
\enqs
As long as  $\tilde\tau_{k}^{(n+1)}$ $\leq$ $\tilde\tau_{1}^{(n+1)}+mh$,  increment $k$ $\rightarrow$ $k+1$~: $\tilde\tau_{k}^{(n+1)}$
$\rightarrow$ $\tilde\tau_{k+1}^{(n+1)}$, until
\beqs
k_{n+1} &=& \sup\big\{ k~:  \tilde\tau_{k}^{(n+1)} \leq \tilde\tau_{1}^{(n+1)}+mh \} \;\; \in \; \{1,\ldots,m\},
\enqs
and continue the induction on $n$~: $n$ $\rightarrow$ $n+1$ until $\tilde\tau_1^{(n+1)} + mh$ $>$ $T$.
\end{itemize}

The optimal impulse control  is given by the finite sequence $\{(\tilde\tau_k^{(n)},\tilde e_k^{(n)})_{k_{n-1}\leq k \leq k_n}, n=0,\ldots,N\}$, where $N$ $=$ $\inf\{ n\geq 0~: \tilde\tau_1^{(n)} +mh> T\}$, and we set by convention $k_{-1}$ $=$ $1$.

\section{Proofs of main results}

\setcounter{equation}{0} \setcounter{Assumption}{0}
\setcounter{Theorem}{0} \setcounter{Proposition}{0}
\setcounter{Corollary}{0} \setcounter{Lemma}{0}
\setcounter{Definition}{0} \setcounter{Remark}{0}

\subsection{Dynamic programming principle}

From the dynamics \reff{eqX} of the controlled process, we  derive easily  the following properties (recall the notations \reff{defl}-\reff{defk}-\reff{defp})~:

\vspace{1mm}

\noindent $\bullet$ Markov property of the pair $(X^\alpha,p(.,\alpha))$ for any $\alpha$ $\in$ $\Ac$, in the sense that
\beqs \label{markov}
\E\Big[\varphi(X_{\theta_2}^\alpha) \Big| \Fc_{\theta_1}\Big] &=& \E\Big[\varphi(X_{\theta_2}^\alpha) \Big| (X_{\theta_1}^\alpha,p(\theta_1,\alpha)) \Big],
\enqs
for any bounded measurable function $\varphi$, and stopping times $\theta_1$ $\leq$ $\theta_2$ a.s.

\vspace{1mm}

\noindent $\bullet$ Causality of the control, in the sense that for any $\alpha$ $=$ $(\tau_i,\xi_i)_{i\geq 1}$ $\in$ $\Ac$, and $\theta$ stopping time,
\beqs \label{causa}
\alpha^\theta & \in & \Ac_{\theta,p(\theta,\alpha)}, \;\;\; \mbox{ and } \;\; p(\theta,\alpha) \; \in \; k(\theta,\alpha) \;\;\; a.s.
\enqs
where we set $\alpha^\theta$ $=$ $(\tau_{i+\iota(\theta,\alpha)},\xi_{i+\iota(\theta,\alpha)})_{i\geq 1}$.

\vspace{1mm}

\noindent $\bullet$ Pathwise uniqueness of the state process,
\beqs \label{pathuni}
X^{t,x,p,\alpha} &=& X^{\theta,X_\theta^{t,x,p,\alpha},p(\theta,\alpha),\alpha^\theta} \;\;\; \mbox{ on } \;\; [\theta,T],
\enqs
for any $(t,x,p)$ $\in$ $\Dc_k$, $k$ $=$ $0,\ldots,m$, $\alpha$ $\in$ $\Ac_{t,p}$, and $\theta$ $\in$ $\Tc_{t,T}$ the set of stopping times   valued in $[t,T]$.

\vspace{2mm}

From the above properties, we deduce  by usual arguments the inequality  \reff{reldynpro2} of the dynamic programming principle, which can be formulated equivalently in

\begin{Proposition} \label{proprogdyn1}
For all $k$ $=$ $0,\ldots,m$, $(t,x,p)$ $\in$ $\Dc_k$, we have
\beqs
v_k(t,x,p) & \leq & \sup_{\alpha\in\Ac_{t,p}} \inf_{\theta\in\Tc_{t,T}}
\E\Big[ \int_t^\theta f(X_s^{t,x,p,\alpha}) ds +  \sum_{t < \tau_i+mh \leq \theta} c(X_{(\tau_i+mh)^-}^{t,x,p,\alpha},\xi_i) \nonumber  \\
& &   \hspace{25mm}  + \;    v_{k(\theta,\alpha)}(\theta,X_\theta^{t,x,p,\alpha},p(\theta,\alpha)) \Big].
\enqs
\end{Proposition}
{\bf Proof.} Fix $(t,x,p)$ $\in$ $\Dc_k$, $k$ $=$ $0,\ldots,m$, and take arbitrary $\alpha$ $\in$ $\Ac_{t,p}$, $\theta$ $\in$ $\Tc_{t,T}$.  
From the definitions of the performance criterion and the value functions, the law of iterated conditional expectations, Markov property, pathwise uniqueness, and causality features of our model, we get the successive relations
\beqs
J_k(t,x,p,\alpha) &=&
\E\Big[ \int_t^\theta f(X_s^{t,x,p,\alpha}) ds +  \sum_{t < \tau_i+mh \leq \theta} c(X_{(\tau_i+mh)^-}^{t,x,p,\alpha},\xi_i)  \\
& & \; + \;  \E\Big[ \int_\theta^T  f(X_s^{t,x,p,\alpha}) ds +  g(X_T^{t,x,p,\alpha})
+  \sum_{\theta < \tau_i+mh \leq T} c(X_{(\tau_i+mh)^-}^{t,x,p,\alpha},\xi_i) \Big| \Fc_\theta \Big] \Big] \\
&=& \E\Big[ \int_t^\theta f(X_s^{t,x,p,\alpha}) ds +  \sum_{t < \tau_i+mh \leq \theta} c(X_{(\tau_i+mh)^-}^{t,x,p,\alpha},\xi_i)  \\
& & \;\;\;\;\;\;\;\;\;  + \; J_{k(\theta,\alpha)}(\theta,X_\theta^{t,x,p,\alpha},p(\theta,\alpha),\alpha^\theta) \Big] \\
& \leq &  \E\Big[ \int_t^\theta f(X_s^{t,x,p,\alpha}) ds +  \sum_{t < \tau_i+mh \leq \theta} c(X_{(\tau_i+mh)^-}^{t,x,p,\alpha},\xi_i)  \\
& & \;\;\;\;\;\;\;\;\;  + \; v_{k(\theta,\alpha)}(\theta,X_\theta^{t,x,p,\alpha},p(\theta,\alpha)) \Big].
\enqs
Since $\theta$ and $\alpha$ are arbitrary, we obtain the required inequality.
\ep

\vspace{1mm}

As usual, the  inequality  \reff{reldynpro1} of the dynamic programming principle requires in addition to the  Markov, causality and pathwise uniqueness properties, a measurable selection theorem. This inequality can be formulated equivalently in

\begin{Proposition}
For all $k$ $=$ $0,\ldots,m$, $(t,x,p)$ $\in$ $\Dc_k$, we have
\beqs
v_k(t,x,p) & \geq & \sup_{\alpha\in\Ac_{t,p}} \sup_{\theta\in\Tc_{t,T}}
\E\Big[ \int_t^\theta f(X_s^{t,x,p,\alpha}) ds +  \sum_{t < \tau_i+mh \leq \theta} c(X_{(\tau_i+mh)^-}^{t,x,p,\alpha},\xi_i) \nonumber  \\
& &   \hspace{25mm}  + \;    v_{k(\theta,\alpha)}(\theta,X_\theta^{t,x,p,\alpha},p(\theta,\alpha)) \Big].
\enqs
\end{Proposition}
{\bf Proof.} Fix $(t,x,p)$ $\in$ $\Dc_k$, $k$ $=$ $0,\ldots,m$,  and arbitrary $\alpha$ $\in$ $\Ac_{t,p}$, $\theta$ $\in$ $\Tc_{t,T}$.
By definition of the value functions, for any $\eps$ $>$ $0$ and $\omega$ $\in$ $\Omega$, there exists
$\alpha_{\eps,\omega}$ $\in$ $\Ac_{\theta(\omega),p(\theta(\omega),\alpha(\omega))}$, which is an $\eps$-optimal control for
$v_{k(\theta(\omega),\alpha(\omega))}$ at  $(\theta,X_\theta^{t,x,p,\alpha},p(\theta,\alpha))(\omega)$. By a measurable selection theorem
(see e.g. Chapter 7 in \cite{bershr78}), there exists $\bar \alpha_\eps$ $\in$ $\Ac_{\theta,p(\theta,\alpha)}$ s.t. $\bar \alpha_\eps(\omega)$
$=$ $\alpha_{\eps,\omega}(\omega)$ a.s., and so
\beq \label{vinfJ}
v_{k(\theta,\alpha)}(\theta,X_\theta^{t,x,p,\alpha},p(\theta,\alpha)) - \eps & \leq &
J_{k(\theta,\alpha)} (\theta,X_\theta^{t,x,p,\alpha},p(\theta,\alpha),\bar\alpha_\eps) \;\;\;\;\;  a.s.
\enq
Now, we define by concatenation the impulse control $\bar\alpha$ consisting of the impulse control components of $\alpha$
until (including eventually) time $\tau$, and the impulse control components of $\bar\alpha_\eps$ strictly after time $\tau$.
By construction, $\bar\alpha$ $\in$ $\Ac_{t,p}$, $X^{t,x,p,\bar\alpha}$ $=$ $X^{t,x,p,\alpha}$ on $[t,\theta]$,
$k(\theta,\bar\alpha)$ $=$ $k(\theta,\alpha)$, $p(\theta,\bar\alpha)$ $=$ $p(\theta,\alpha)$, and $\bar\alpha^\theta$ $=$ $\bar\alpha_\eps$.
Hence, similarly as in Proposition \ref{proprogdyn1}, by using law of iterated conditional expectations, Markov property, pathwise uniqueness, and causality features of our model, we get
\beqs
J_k(t,x,p,\bar\alpha) &=& \E\Big[ \int_t^\theta f(X_s^{t,x,p,\alpha}) ds +  \sum_{t < \tau_i+mh \leq \theta} c(X_{(\tau_i+mh)^-}^{t,x,p,\alpha},\xi_i)  \\
& & \;\;\;\;\;\;\;\;\;  + \; J_{k(\theta,\alpha)}(\theta,X_\theta^{t,x,p,\alpha},p(\theta,\alpha),\bar\alpha_\eps) \Big].
\enqs
Together with \reff{vinfJ}, this implies
\beqs
v_k(t,x,p) \; \geq \; J_k(t,x,p,\bar\alpha) & \geq & \E\Big[ \int_t^\theta f(X_s^{t,x,p,\alpha}) ds +  \sum_{t < \tau_i+mh \leq \theta} c(X_{(\tau_i+mh)^-}^{t,x,p,\alpha},\xi_i)  \\
& & \;\;\;\;\;\;\;\;\;  + \;  v_{k(\theta,\alpha)}(\theta,X_\theta^{t,x,p,\alpha},p(\theta,\alpha)) \Big] - \eps.
\enqs
From the arbitrariness of $\eps$, $\alpha$, and $\theta$, this proves the required result.
\ep

\vspace{1mm}

We end this paragraph by proving Corollary \ref{coroldyn}.

\vspace{1mm}

\noindent {\bf Proof of Corollary \ref{coroldyn}.}

\noindent (i) Fix  $k$ $\in$ $\{1,\ldots,m\}$,  $(t,x)$ $\in$ $[0,T]\times\R^d$,
$p$  $=$ $(t_i,e_i)_{1\leq i\leq k}$ $\in$ $P_t^1(k)$ such that  $t_1+mh$ $\leq$ $T$, and $\theta$
stopping time valued in  $[t,(t_k+h)\wedge(t_1+mh))$. Then, we observe that for all $\alpha$ $=$ $(\tau_i,\xi_i)_{i\geq 1}$  $\in$ $\Ac_{t,p}$, $X^{t,x,p,\alpha}$ $=$ $X^{t,x,0}$ on $[t,\theta]$, $\tau_i+mh$ $>$ $\theta$, $k(\theta,\alpha)$ $=$ $k$, and $p(\theta,\alpha)$ $=$ $p$ a.s. Hence, relation \reff{dynprovkpc1} follows  immediately from \reff{reldynpro}.

\vspace{1mm}

\noindent (ii)   For  $k$  $\in$ $\{0,\ldots,m-1\}$, $p$ $=$ $(t_i,e_i)_{1\leq i\leq k}$  $\in$ $P_t^2(k)$ such that  $t_1+mh$ $\leq$ $T$,
and $\theta$ stopping time valued in $[t,(t_1+mh)\wedge (t+h))$. Let $\alpha$ $=$ $(\tau_i,\xi_i)_{i\geq 1}$  be some arbitrary element in $\Ac_{t,p}$,
and set $\tau$ $=$ $\tau_{k+1}$, $\xi$ $=$ $\xi_{k+1}$. Notice that $(\tau,\xi)$ $\in$ $\Ic_t$.
Then, we see that
$X^{t,x,p,\alpha}$ $=$ $X^{t,x,0}$ on $[t,\theta]$, $\tau_i+mh$ $>$ $\theta$, $k(\theta,\alpha)$ $=$ $k$, $p(\theta,\alpha)$ $=$ $p$ if $\theta$ $<$ $\tau$, and  $k(\theta,\alpha)$ $=$ $k+1$, $p(\theta,\alpha)$ $=$ $p\cup (\tau,\xi)$ if $\theta$ $\geq$ $\tau$.
We deduce from \reff{reldynpro1} that
\beqs
v_k(t,x,p) & \geq &  \E\Big[ \int_t^{\theta} f(X_s^{t,x,0}) ds +
v_k(\theta,X_{\theta}^{t,x,0},p)  1_{\theta < \tau}   \\
& &   \hspace{25mm} + \; v_{k+1}(\theta,X_{\theta}^{t,x,0}, p \cup (\tau,\xi))   1_{\tau \leq \theta} \Big],  \nonumber
\enqs
and this inequality holds for any $(\tau,\xi)$ $\in$ $\Ic_t$ by arbitrariness of $\alpha$. Furthermore, from \reff{reldynpro2},
for all $\eps$ $>$ $0$, there exists $(\tau,\xi)$ $\in$ $\Ic_t$ s.t.
\beqs
v_k(t,x,p) -\eps & \leq &  \E\Big[ \int_t^{\theta} f(X_s^{t,x,0}) ds +
v_k(\theta,X_{\theta}^{t,x,0},p)  1_{\theta < \tau}   \\
& &   \hspace{25mm} + \; v_{k+1}(\theta,X_{\theta}^{t,x,0}, p \cup (\tau,\xi))   1_{\tau \leq \theta} \Big].   \nonumber
\enqs
The two previous inequalities give the required relation
\beqs
v_k(t,x,p) & = &  \sup_{(\tau,\xi)\in\Ic_t} \E\Big[ \int_t^{\theta} f(X_s^{t,x,0}) ds +
v_k(\theta,X_{\theta}^{t,x,0},p)  1_{\theta < \tau}  \label{dynprovkpc22} \\
& &   \hspace{25mm} + \; v_{k+1}(\theta,X_{\theta}^{t,x,0}, p \cup (\tau,\xi))   1_{\tau \leq \theta} \Big].
\enqs

\subsection{Viscosity properties} \label{paravisco}

In this paragraph,   we  prove the viscosity property stated in  Proposition \ref{provisco}.  
We first  state an auxiliary result, which can be proved  similarly  as in Lemma 5.1 in  \cite{lyvmnipha07}.  For any
locally bounded function $u$ on $\Dc_{k+1}^{m}$,  $k$ $=$ $0,\ldots,m-1$, we define the  locally bounded function  $\Hc u$ on $\Dc_{k}^{2,m}$ by
$\Hc u (t,x,p)$ $=$ $\sup_{e\in E} u(t,x,p\cup (t,e))$.

\begin{Lemma} \label{lemmeVathana}
Let $u$ be a locally bounded function on $\Dc_{k+1}^{m}$,  $k$ $=$ $0,\ldots,m-1$. Then,  $\Hc \overline{u}$ is upper-semicontinuous, 
and  $\overline{\Hc u}$ $\leq$ $\Hc \overline{u}$.
\end{Lemma}
{\bf Proof.}
Fix some $(t,x,p)$ $\in$ $\Dc_{k}^{2,m}$, and let $(t_n,x_n,p_n)_{n\geq 1}$ be a sequence in $\Dc_{k}^{2,m}$ converging to $(t,x,p)$ as $n$ goes to infinity.  Since $\overline{u}$ is upper-semicontinuous, and $E$ is compact,  there exists a sequence $(e_n)_n$ valued in $E$, such that 
\beqs
\Hc \overline{u}(t_n,x_n,p_n) &=& \overline{u}(t_n,x_n,p_n \cup (t_n,e_n)), \;\;\; n \geq 1. 
\enqs
The sequence $(e_n)_n$ converges, up to a subsequence, to some $\hat e$ $\in$ $E$, and so
\beqs
\Hc \overline{u}(t,x,p) \; \geq \; \overline{u}(t,x,p\cup (t,\hat e)) & \geq & 
\limsup_{n\rightarrow\infty} \overline{u}(t_n,x_n,p_n\cup (t_n,e_n)) \; = \;  \limsup_{n\rightarrow\infty} \Hc \overline{u}(t_n,x_n,p_n), 
\enqs
which shows that $\Hc\overline{u}$ is upper-semicontinuous. 

On the other hand,  fix some $(t,x,p)$ $\in$ $\Dc_{k}^{2,m}$, and let $(t_n,x_n,p_n)_{n\geq 1}$ be a sequence in $\Dc_{k}^{2,m}$ converging to $(t,x,p)$ 
s.t. $\Hc u(t_n,x_n,p_n)$ converges to  $\overline{\Hc u}(t,x,p)$.  Then, we have 
\beqs
\overline{\Hc u}(t,x,p) \; = \; \lim_{n\rightarrow\infty} \Hc u(t_n,x_n,p_n) & \leq & \limsup_{n\rightarrow\infty} \Hc \overline{u}(t_n,x_n,p_n) \; \leq \; 
\Hc \overline{u}(t,x,p),
\enqs
which shows that $\overline{\Hc u}$ $\leq$ $\Hc \overline{u}$.
\ep

\vspace{2mm}

Now, we prove the sub and supersolution property of the family  $v_{k}$, $k$ $=$ $0,\ldots,m$.  There is no difficulty on the domain
$\Dc_{k}^{1,m}$ since  locally  no  impulse control is possible. Hence, in this case,
the viscosity properties can be derived as for an uncontrolled state process, and the proof is standard from the dynamic programming principle 
\reff{dynprovkpc1}, see  e.g. \cite{pham05}.  Notice  that  since  the domain $\T_p^{1}(k)$ is of the form $[t_{k},(t_{k}+h)\wedge(t_{1}+mh))$, we have no problem at the boundary. Indeed, this set is open at $(t_{k}+h)\wedge(t_{1}+mh)$, which is the usual situation, and the closedness at $t_{k}$ does not introduce difficulties, as the value function is not defined before $t_{k}$.  Hence, when taking approximations of the upper and lower semicontinous envelopes of $v_{k}$, we  only  need to consider points of the domain such that $t\geq t_{k}$, where the dynamic programming relation  
\reff{dynprovkpc1} holds. The  proof of the viscosity property of the value functions $v_k$ to  \reff{viscoQ2} on $\Dc_{k}^{2,m}$ is more subtle.  Indeed, in addition to the specific form of equation  \reff{viscoQ2}, we have to carefully address  the discontinuity of the PDE system  \reff{viscoQ1}-\reff{viscoQ2} on the left boundary of $\T_p^{2}(k)$. In the sequel, we  focus on the domain  $\Dc_{k}^{2,m}$, $k$ $=$ $0,\ldots,m-1$.

\vspace{2mm}

\noindent {\bf Proof of the supersolution property on} $\Dc_{k}^{2,m}$.

\vspace{1mm}

\noindent We first prove that  for $k$ $=$ $0,\ldots,m-1$,  $(t_0,x_0,p_0)\in \Dc_{k}^{2,m}$~:
\beq \label{aProuverSursol1}
\underline{v_{k}}(t_0,x_0,p_0)\geq \sup_{e \in E}\underline{v_{k+1}}(t_0,x_0,p_0\cup(t_0,e)).
\enq
By definition of $\underline{v_{k}}$, there exists a sequence  $(t_{n},x_{n},p_{n})_{n\geq 1} \in \Dc_{k}^{m}$ such that~:
\beq \label{deflsc}
v_{k}(t_{n},x_{n},p_{n}) \to \underline{v_{k}}(t_0,x_0,p_0) \;\; \text{ with } \;\;  (t_{n},x_{n},p_{n}) \to (t_0,x_0,p_0).
\enq
We set  $p_0$ $=$  $(t_i^0,e_i^0)_{1\leq i\leq k}$,  $p_n$ $=$ $(t_i^n,e_i^n)_{1\leq i\leq k}$, and we distinguish the two following cases~:
\begin{itemize}
\item  If $t_0$ $>$ $t_{k}^0+h$, then, for $n$ sufficiently large, we have $t_{n} \geq t_{k}^{n}+h$, i.e. $p_n$ $\in$ $P_{t_n}^2(k)$. Hence, 
from the dynamic programming principle by making an immediate impulse control, i.e. by applying  \reff{dynprovkpc2} to  
$v_k(t_n,x_n,p_n)$ with   $\theta$ $=$  $\tau$ $=$ $t_n$, and $e$ $\in$ $E$,  we have
\beqs
v_{k}(t_{n},x_{n},p_{n}) &\geq&  v_{k+1}(t_{n},x_{n},p_{n}\cup(t_{n},e)) \; \geq \;  \underline{v_{k+1}}(t_{n},x_{n},p_{n}\cup(t_{n},e)).  
\enqs
By sending $n$ to infinity with \reff{deflsc}, and  since   $\underline{v_{k+1}}$ is lower-semicontinuous,  we obtain the required relation 
\reff{aProuverSursol1} from the arbitrariness of $e$ in $E$.  

\item  if $t_0$ $=$ $t_{k}^0+h$,  we apply the dynamic programming principle by making an impulse control as soon as possible. This means that in relation \reff{reldynpro1} for $v_k(t_n,x_n,p_n)$, we choose $\alpha$ $=$ 
$(\tau_i,\xi_i)_{i\geq 1}$ $\in$ $\Ac_{t_n,p_n}$, $\theta$ $=$ $\tau_{k+1}$ $=$ $\theta_n$ $:=$  $t_{n} \vee (t_{k}^{n}+h)$, $\xi_{k+1}$ $=$ $e$ $\in$ $E$, so that~: 
\beqs
v_{k}(t_{n},x_{n},p_{n}) &\geq&  \E\Big[  \int_{t_{n}}^{ \theta_n}f(X_{s}^{n})ds 
+ \sum_{t_n < \tau_i+mh\leq \theta_n} c(X_{(\tau_i+mh)^-}^n,\xi_i)  \\
& &  \;\;\;\;\;\;\;  +  \;  \underline{v_{k+1}}(\theta_n,X_{\theta_n}^{n},p_n\cup (\theta_n,e)) \big].
\enqs
Here $X^n$ $:=$ $X^{t_n,x_n,0}$.  Since $t_n$, $\theta_n$ $\to$ $t_0$, $p_n$ $\to$ $p_0$, $X_{\theta_n}^{n} \to x_0$ a.s.,  as $n$ goes to infinity,  and from estimate \reff{estimX2} and  the linear growth condition on $f$, $c$, $\underline{v_{k+1}}$,  
we can use the dominated convergence theorem to obtain~: 
\beqs
\underline{v_{k}}(t_0,x_0,p_0)\geq \underline{v_{k+1}}(t_0 ,x_0,p_0\cup(t_0,e)),
\enqs
which implies \reff{aProuverSursol1} from the arbitrariness of $e$ $\in$ $E$.
\end{itemize}
Finaly, in order  to complete  the viscosity   supersolution property of $v_k$ to  \reff{viscoQ2} on $\Dc_{k}^{2,m}$,  it remains  to
show that $v_k$ is  a supersolution to~:
\beqs
- \Dt{v_k}(t,x,p)  - \Lc v_k(t,x,p) - f(x) &\geq& 0,
\enqs
on $\Dc_{k}^{2,m}$. This proof is standard by  using the dynamic programming relation  \reff{dynprovkpc2} with $\tau$ $=$ $\infty$ and It\^o's formula,
see   \cite{pham05} for the details.
\ep

\vspace{3mm}

\noindent {\bf Proof of the subsolution property on } $\Dc_{k}^{2,m}$.

\vspace{1mm}

\noindent We follow arguments in  \cite{lyvmnipha07}.
Let $(t_0,x_0,p_0)$ $\in$ $\Dc_{k}^{2,m}$ and $\varphi$ $\in$ $C^{1,2}(\Dc_{k}^{2,m})$ such that
$\overline{v_{k}}(t_0,x_0,p_0)$ $=$ $\varphi(t_0,x_0,p_0)$ and $\varphi\geq\overline{v_{k}}$ on $\Dc_{k}^{2,m}$.  If
$\overline{v_{k}}(t_0,x_0,p_0)$  $\leq$ $\Hc \overline{v_{k+1}}(t_0,x_0,p_0)$,  then the subsolution inequality holds trivially. Now, if
$\overline{v_{k}}(t_0,x_0,p_0)$ $>$ $\Hc \overline{v_{k+1}}(t_0,x_0,p_0)$,  we argue by contradiction by assuming on the contrary that
\beqs
\eta \; := \;  - \Dt{\varphi}(t_0,x_0,p_0) - \Lc \varphi(t_0,x_0,p_0)- f(x_0) & > & 0.
\enqs
We set  $p_0$ $=$  $(t_i^0,e_i^0)_{1\leq i\leq k}$.  By continuity of $\varphi$ and its derivatives, there exists some $\delta>0$ with 
$t_0+\delta$ $<$ $(t_1^0+mh)\wedge T$ such that~:
\beq \label{varphieta}
- \Dt{\varphi}   - \Lc \varphi  - f  &>& \frac{\eta}{2},  \;\;\; \mbox{ on }  \;
((t_0-\delta,t_0+\delta)\times B(x_0,\delta)\times B(p_0,\delta)) \cap \Dc_{k}^{2,m}.
\enq
From the definition of $\overline{v_{k}}$, there exists a sequence $(t_{n},x_{n},p_{n})_{n\geq 1 }$ $\in$
$((t_0-\delta,t_0+\delta)\times B(x_0,\delta)\times B(p_0,\delta))$ $\cap$ $\Dc_{k}^{2,m}$
such that $(t_{n},x_{n},p_{n}) \to (t_0,x_0,p_0)$ and $v_{k}(t_{n},x_{n},p_{n}) \to \overline{v_{k}}(t_0,x_0,p_0)$ as $n\to  \infty$. By continuity of
$\varphi$ we also have that $\gamma_{n}$ $:=$ $v_{k}(t_{n},x_{n},p_{n}) - \varphi(t_{n},x_{n},p_{n})$ converges to $0$ as $n \to \infty$.
We set $p_n$ $=$ $(t_i^n,e_i^n)_{1\leq i\leq k}$.
From  the dynamic programming principle \reff{dynprovkpc2}, for each $n$ $\geq$ $1$, there exists a control
$(\tau^{n},\xi^{n})$ $\in$ $\Ic_{t_n}$ such that
\beq
v_{k}(t_{n},x_{n},p_{n})  -  \frac{\eta}{4} \delta_n &\leq& \E \left[\int_{t_{n}}^{\theta_n}f(X_{s}^{n})ds+ v_{k}(\theta_n,X_{\theta_n}^{n},p_{n}) 1_{\theta_n <\tau_{n}} \right.
\nonumber \\
& &  \left. \;\;\;\;\; + \;  v_{k+1}(\theta_n,X_{\theta_n},p_{n}\cup(\tau_{n},\xi_{n}))1_{\tau_{n}\leq\theta_n} \right] . \label{uneInegalite}
\enq
Here $X^{n}$ $:=$ $X^{t_{n},x_{n},0}$,  we  choose  $\theta_{n}$ $=$ $\vartheta_n \wedge (t_n+\delta_n)$, with  
$\vartheta_n$ $=$ $ \inf\{s\geq t_{n} :  X_{s}^{n} \notin B(x_{n},\frac{\delta}{2})\}$,  and $(\delta_n)_n$ is a strictly positive sequence such that 
\beqs
\delta_n \rightarrow 0, & & \frac{\gamma_n}{\delta_n} \; \rightarrow \; 0, \;\;\; \mbox{ as } \; n \rightarrow \infty.
\enqs
On the other hand,  from  Lemma \ref{lemmeVathana}, we have
\beqs
\overline{\Hc v_{k+1}}(t_0,x_0,p_0)  &\leq&  \Hc \overline{v_{k+1}} (t_0,x_0,p_0) \; < \; \overline{v_{k}}(t_0,x_0,p_0)
 \; \leq \;  \varphi(t_0,x_0,p_0).
\enqs
Hence, since $\overline{\Hc v_{k+1}}$  is u.s.c. and $\varphi$ is continuous, the inequality $\Hc v_{k+1}$ $\leq$ $\varphi$  holds in a neighborhood of $(t_0,x_0,p_0)$, and so  for  sufficiently large $n$, we get~:
\beqs
v_{k+1}(\theta_{n},X_{\theta_{n}}^{n},p_{n}\cup(\tau_{n},\xi_{n} ))1_{\tau_{n}\leq
\theta_{n}}  &\leq&   \varphi(\theta_{n},X_{\theta_{n}}^{n},p_{n})1_{\tau_{n} \leq \theta_{n}}  \;\;\; a.s. 
\enqs
Together with \reff{uneInegalite}, this yields~:
\beqs
\varphi(t_{n},x_{n},p_{n}) \ +\gamma_{n} - \frac{\eta}{4} \delta_n &\leq&
\E \left[\int_{t_{n}}^{\theta_{n}}f(X_{s}^{n})ds  +  \varphi(\theta_{n},X_{\theta_{n}}^{n},p_{n} )\right].
\enqs
By applying It\^o's formula to $\varphi(s,X_{s}^{n},p_{n})$ between $s$ $=$ $t_{n}$ and $s$ $=$ $\theta_{n}$, and dividing by $\delta_n$, we then get~:
\beq
\frac{\gamma_n}{\delta_n}  - \frac{\eta}{4} &\leq&  \frac{1}{\delta_n} \E \left[\int_{t_{n}}^{\theta_{n}}\left(\Dt{\varphi}+\Lc \varphi +f \right)(s,X_{s}^{n},p_{n})ds \right] \; \leq \;   -\frac{\eta}{2}\E\left[\frac{\theta_{n}-t_{n}}{\delta_n} \right], \label{itophi}
\enq
from \reff{varphieta}.  Now, from the growth linear condition on $b$, $\sigma$, Burkholder-Davis-Gundy inequality and Gronwall's lemma, we have  the standard estimate~: $\E[\sup_{s\in [t_n,t_n+\delta_n]}|X_s^n-x_n|^2]$ $\to$ $0$, so that by Chebichev inequality, $\P[\vartheta_n\leq t_n+\delta_n]$ $\to$ 
$0$, as $n$ goes to infinity,  and therefore by definition of $\theta_n$~: 
\beqs
1 \; \geq \;  \E\left[\frac{\theta_{n}-t_{n}}{\delta_n} \right] & \geq & \P[\vartheta_n > t_n+\delta_n] \; \rightarrow \; 1,  \;\; \mbox{ as } \; n \to \infty.
\enqs
By sending $n$  to infinity into \reff{itophi},  we obtain the required contradiction~: $-\frac{\eta}{4}$ $\leq$ $-\frac{\eta}{2}$. 
\ep

\subsection{Sequential comparison results} \label{paracomp}

In this paragraph, we  prove sequential comparison results. We
consider the  sets $\Theta_k(n)$,  $\T^n(0)$, $\Dc_k(n)$, $\Dc_k^m(n)$, $\Dc_k^{i,m}(n)$,
introduced in Section \ref{secalgo} for $k$ $=$ $0,\ldots,m$, and  $n$ $=$ $m,\ldots,M$, and  we define sequential viscosity
solutions  as follows.

\begin{Definition} Let $n$ $\in$ $\{m+1,\ldots,N\}$.
We say that  a family of locally bounded  functions $w_k$ on $\Dc_k^m(n)$, $k$ $=$ $0,\ldots,m$,  is a viscosity supersolution  (resp. subsolution)  of
\reff{viscoQ1}-\reff{viscoQ2} at step $n$ if~:

\noindent (i)  for all $k$ $=$ $1,\ldots,m$,   $(t_0,x_0,p_0)$ $\in$ $\Dc_k^{1,m}(n)$, and $\varphi$ $\in$
$C^{1,2}(\Dc_k^{1,m}(n))$, which realizes a local minimum of $\underline{w_k}-\varphi$  (resp. maximum of
$\overline{w_k}-\varphi$), we have
\beqs
 - \Dt{\varphi}(t_0,x_0,p_0)  - \Lc \varphi(t_0,x_0,p_0)  - f(x_0)   & \geq  & 0  \;\; (resp. \; \leq \; 0).
\enqs
(ii) for all $k$ $=$ $0,\ldots,m-1$,    $(t_0,x_0,p_0)$ $\in$ $\Dc_k^{2,m}(n)$,  and $\varphi$ $\in$
$C^{1,2}(\Dc_k^{2,m}(n))$, which realizes a local minimum of $\underline{w_k}-\varphi$  (resp. maximum of
$\overline{w_k}-\varphi$), we have
\beqs
 \min\big\{  - \Dt{\varphi}(t_0,x_0,p_0)  - \Lc \varphi(t_0,x_0,p_0)  - f(x_0) \; , \; \hspace{7mm}  & &   \\
\;\;\;\;\;\;\;   \underline{w_k}(t_0,x_0,p_0) -  \sup_{e\in E} \underline{w_{k+1}}(t_0,x_0,p_0 \cup (t_0,e)) \big\} &  \geq &   0
\enqs
(resp.
\beqs
 \min\big\{  - \Dt{\varphi}(t_0,x_0,p_0)  - \Lc \varphi(t_0,x_0,p_0)  - f(x_0) \; , \; \hspace{7mm}  & &   \\
\;\;\;\;\;\;\;   \overline{w_k}(t_0,x_0,p_0) -  \sup_{e\in E} \overline{w_{k+1}}(t_0,x_0,p_0 \cup (t_0,e)) \big\} &  \leq &   0 ).
\enqs
We say that a family of locally bounded  functions $w_k$ on $\Dc_k^m(n)$, $k$ $=$ $0,\ldots,m$,  is a
viscosity solution of  \reff{viscoQ1}-\reff{viscoQ2} at step
$n$ if it is a viscosity supersolution and subsolution of  \reff{viscoQ1}-\reff{viscoQ2} at step $n$.
\end{Definition}

We then prove the following comparison principle at step $n$.

\begin{Proposition} \label{procompn}  Let $n$ $\in$ $\{m+1,\ldots,N\}$.
Let $u_k$   (resp. $w_k$), $k$ $=$ $0,\ldots,m$, be a family of  viscosity subsolution (resp. supersolution) of \reff{viscoQ1}-\reff{viscoQ2} at step $n$
satisfying growth condition \reff{growthvk}. Suppose also that $w_{k}$ satisfies \reff{condSupplDisc}.
If $u_{k}$ and $w_{k}$ are such that for all $x$ $\in$ $\R^d$
\beqs
\overline{u_k}(t_1+mh,x,p) & \leq & \underline{w_k}(t_1+mh,x,p), \;\; \; p = (t_i,e_i)_{1\leq i\leq k}  \in \Theta_k^m(n)\times E^k, \; k \geq  1, \\
 \overline{u_0}(T,x)  & \leq &   \underline{w_0}(T,x).
\enqs
Then,  $\overline{u_k}$ $\leq$ $\underline{w_k}$ on $\Dc_k^m(n)$,  $k$ $=$ $0,\ldots,m$.
\end{Proposition}

\begin{Remark}
{\rm   We recall some  basic  definitions and properties  in viscosity solutions theory, which shall be used in the proof of the above proposition.  
Consider the general PDE 
\beq \label{PDEF}
F(t,x,w,\Dt{w},D_x w,D_x^2 w) &=& 0 \;\;\; \mbox{ on } \;  [t_0,t_1)\times \Oc,
\enq
where $t_0$ $<$ $t_1$, and $\Oc$ is an open set in $\R^d$.  There is an equivalent definition of viscosity solutions to \reff{PDEF} in terms of semi-jets 
$\bar J^{2,+} w(t,x)$  and $\bar J^{2,-} w(t,x)$ associated respectively to an upper-semicontinuous (u.s.c.) and lower-semicontinuous (l.s.c.) 
function $w$ (see \cite{craishlio92} or \cite{fleson93} for the definition of semi-jets)~:  
an u.s.c. (resp. l.s.c.) function $w$ is a viscosity subsolution (resp. supersolution) to \reff{PDEF} if and only if for all $(t,x)$ $\in$ $[t_0,t_1)\times\Oc$, 
\beqs
F(t,x,w(t,x),r,q,A) & \leq \; ( \mbox{ resp. } \geq) & 0,  \;\;\; \forall (r,q,A) \in  \bar J^{2,+} w(t,x) \; (\mbox{ resp. } \;   \bar J^{2,-} w(t,x)).  
\enqs  
For $\eta$ $>$ $0$, we say that  $w^{\eta}$ is a viscosity  $\eta$-strict supersolution to  \reff{PDEF},  if $w^\eta$ is a viscosity supersolution to
\beqs
F(t,x,w^{\eta},\Dt{w^\eta},D_xw^{\eta},D_x^{2}w^{\eta})  & \geq & \eta,  \;\;\; \mbox{ on } \;  [t_0,t_1)\times \Oc.
\enqs
in the sense that it is a viscosity supersolution to $F(t,x,w^{\eta},\Dt{w^\eta},D_xw^{\eta},D_x^{2}w^{\eta})$ $-$ $\eta$ $=$ $0$, 
on  $[t_0,t_1)\times \Oc$. 
}
\end{Remark}

\vspace{2mm}

As usual when dealing with variational inequalities,   we begin the proof of the
comparison principle by showing the existence of viscosity $\eta$-strict supersolutions for equation \reff{viscoQ1}-\reff{viscoQ2}.

\begin{Lemma} \label{existStrict}
Let $w_{k}$, $k$ $=$ $0,\ldots,m$,  be a family of viscosity supersolutions of \reff{viscoQ1}-\reff{viscoQ2} satisfying \reff{condSupplDisc}.
Then, for any $\eta>0$, there exists a  family of viscosity $\eta$-strict supersolutions $w_{k}^{\eta}$ of \reff{viscoQ1}-\reff{viscoQ2} such that for $k$ $=$
$0,\ldots,m$~:
\beq \label{encadr}
w_k(t,x,p) + \eta C_{1}|x|^{2} & \leq & w_k^\eta(t,x,p) \; \leq \; w_k(t,x,p) + \eta C_2 (1+|x|^2),  \; (t,x,p) \in \Dc_k^m,
\enq
for some positive constants $C_{1}$, $C_{2}$ independent on $\eta$. Moreover, for   $k$ $=$ $0,\ldots,m-1$,
$(t,x,p)\in \Dc_{k}^{m}$, $p$ $=$ $(t_i,e_i)_{1\leq i\leq k}$ with $t=t_{k}+h$, we have~:
\beq \label{condSupplDiscEta}
\underline{w_{k}^{\eta}}(t,x,p)&\geq& \sup_{e\in E} \underline{w_{k+1}^{\eta}}(t,x,p\cup(t,e)) +\eta.
\enq
\end{Lemma}
{\bf Proof.}
For $\eta>0$, consider the functions~:
\beqs
w_{k}^{\eta}(t,x,p) = w_{k}(t,x,p)+\eta\phi_{1,k}(t)+\eta\phi_{2}(t,x), \;\;\;
\phi_{1,k}(t)=\left[(T-t)+(m-k)\right],\\
\phi_{2}(t,x)=\frac{1}{2} e^{L\left(T-t\right)}\left(1 + |x|^{2}  \right),
\enqs
with $L$ a positive constant to be determined later. It is clear that $w_k^\eta$ satisfies \reff{encadr} with $C_1$ $=$ $1/2$ and
$C_2$ $=$ $T+m+e^{LT}/2$. Moreover, we easily show that $w_{k}+\eta\phi_{1,k}^{\eta}$ is a viscosity supersolution to
\beq \label{viscosurphi1}
- \Dt{(w_{k}+\eta\phi_{1,k})}  - \Lc (w_k +\eta\phi_{1,k}) - f  &\geq& \eta.
\enq
This is derived from the fact that $-\Dt{\phi_{1,k}}$ $-$ $\Lc\phi_{1,k}$ $=$ $1$, and $w_{k}$ is a viscosity 
supersolution to $-\Dt{w_{k}}$ $-$ $\Lc w_{k}$ $-$ $f$ $\geq$ $0$.
We now show that $\phi_{2}$ is a supersolution to
\beq \label{viscosurphi2}
- \Dt{ \phi_{2}}  - \Lc\phi_{2}  &\geq& 0.
\enq
This is done by calculating this quantity explicitely. Indeed, we have
\beqs
\Dt{ \phi_{2}}(t,x) = - \frac{L}{2} e^{L(T-t)} (1+|x|^{2}), \;\;\;
\Lc\phi_{2}(t,x)=e^{L(T-t)}\left( b(x).x + {\rm tr}\left(\sigma\sigma'(x) \right)\right).
\enqs
Since $b$ and $\sigma$ are of linear growth, we thus obtain~:
\beqs
- \Dt{ \phi_{2}}(t,x)  - \Lc\phi_{2}(t,x) &\geq&  e^{L(T-t)} \left[\frac{L}{2}(1+|x|^{2})- C(1+|x|+|x|^{2})  \right],
\enqs
for some constant $C$ independent of $t,x$. Therefore, by taking $L$ sufficiently large, we get the required inequality \reff{viscosurphi2},
which shows together with \reff{viscosurphi1} that $w_k^\eta$ is a viscosity supersolution to
\beq \label{wketasur}
- \Dt{w_{k}^\eta}  - \Lc w_k^\eta  - f  &\geq& \eta.
\enq
Moreover, since
\beqs
\underline{w_{k}}(t,x,p)  -   \sup_{e \in E}  \underline{w_{k+1}}(t,x,p\cup (t,e)) &\geq& 0,
\enqs
we immediately get
\beqs
& & \underline{w_{k}^{\eta}}(t,x,p)  -   \sup_{e \in E}  \underline{w_{k+1}^{\eta}}(t,x,p\cup (t,e)) \\
& = & \underline{w_{k}}(t,x,p)+\eta\phi_{1,k}(t)- \sup_{e \in E}  \underline{w_{k+1}}(t,x,p\cup (t,e))-\eta\phi_{1,k+1}(t) \\
& \geq & \eta \phi_{1,k}(t) - \eta\phi_{1,k+1}(t) \;  \geq \;  \eta.
\enqs
Together with \reff{wketasur}, this proves the required viscosity $\eta$-strict supersolution property for $w_{k}^\eta$ to \reff{viscoQ1}-\reff{viscoQ2}.
\ep

\vspace{3mm}

The main step in the  proof of  Proposition \ref{procompn} consists  in the comparison principle for $\eta$-strict supersolutions. Notice from
\reff{encadr} that once $w_k$ satisfies a linear growth condition, then $w_k^\eta$ satisfies the quadratic growth lower-bound condition~:
\beq \label{quadGrowth}
\eta C_{1} \left|x\right|^{2} - C_2 & \leq &  w_{k}^{\eta}(t,x,p), \;\;\; (t,x,p) \in \Dc_k^m,
\enq
for some positive constants $C_1$, $C_2$.

\begin{Lemma} \label{comparStict}  Let $n$ $\in$ $\{m+1,\ldots,N\}$ and $\eta$ $>$ $0$.
Let $u_k$   (resp. $w_k$), $k$ $=$ $0,\ldots,m$, be a family of  viscosity subsolution (resp. $\eta$-strict supersolution)
of \reff{viscoQ1}-\reff{viscoQ2} at step $n$, with $u_{k}$ satisfying the linear growth condition \reff{growthvk} and $w_{k}$ satisfying
the quadratic growth condition
\reff{quadGrowth}. Suppose that  for all $x$ $\in$ $\R^d$,
\beq \label{condtermcomp1}
\overline{u_k}(t_1+mh,x,p) & \leq & \underline{w_k}(t_1+mh,x,p), \;\; p = (t_i,e_i)_{1\leq i\leq k}  \in \Theta_k^m(n)\times E^k,  k \geq  1, \\
 \overline{u_0}(T,x)  & \leq &   \underline{w_0}(T,x).  \label{condtermcomp2}\\
\underline{w}_{k}(t_{k}+h,x,\pi) & \geq & \sup_{e\in E} \underline{w_{k+1}}\left(t_{k}+h,x,p\cup(t_{k}+h,e)\right)+\eta ,\label{condSuppl} \\
  &  &  \;\;\;  p = (t_i,e_i)_{1\leq i\leq k}  \in  \Theta_k^m(n)\times E^k, \; k \leq  m-1. \nonumber
\enq
Then,  $\overline{u_k}$ $\leq$ $\underline{w_k}$ on $\Dc_k^m(n)$,  $k$ $=$ $0,\ldots,m$.
\end{Lemma}
{\bf Proof.}
From the linear growth of $u_k$, and from the quadratic growth lower-bound  of $w_k$, we have
\beqs
\overline{u_{k}}(t,x,p)-\underline{w_{k}}(t,x,p) &\leq& C_{1}\left(1+ \left|x \right| \right)-C_{2}\left|x\right|^{2}, \;\;\; k=0,\ldots,m, \; (t,x,p) \in \Dc_k^m(n),
\enqs
for some positive constants $C_1$, $C_{2}$.  Thus, for all $k$, the supremum of the u.s.c function $\overline{u_k}-\underline{w_k}$ is attained on a compact set  that only depends on $C_{1}$ and $C_{2}$.  
Hence, one can find $k_0$ $\in$ $\{0,\ldots,m\}$, $(t_0,x_0,p_0)$ $\in$ $\Dc_{k_0}^m(n)$ such that~:
\beq 
M & := & \sup_{\tiny\begin{array}{l} k \in \{ 0,\ldots,m \}  \\ (t,x,p) \in \Dc_k^m(n) \end{array}}   
\left[ \overline{u_k}(t,x,p) - \underline{w_k}(t,x,p)\right] \nonumber \\
&=& \overline{u_{k_0}}(t_0,x_0,p_0) - \underline{w_{k_0}} (t_0,x_0,p_0), \label{defM}
\enq
and we have to show that  $M$ $\leq$ $0$.  We set $p_0$ $=$ $(t_i^0,e_i^0)_{1\leq i\leq k_0}$, and we
distinguish the five possible cases concerning $(k_0,t_0,x_0,p_0)$~:
\begin{itemize}
\item \it{Case 1} :  $k_0$ $\neq$ $0$, $t_0$ $=$ $t_{1}^0+mh$.
\item \it{Case 2} : $k_0$ $=$ $0$, $t_0$ $=$ $T$.
\item \it{Case 3} :  $k_0$ $\neq$ $0$,  $t_0$ $\in$ $\T_{p_0}^1(k_0)$.
\item \it{Case 4} :  $k_0$ $=$ $0$, $t_0$ $\in$ $[0,T)$ or $k_0$ $\in$ $\{1,\ldots,m-1\}$, $t_0$ $\in$ $\T_{p_0}^2(k_0)$, $t_0$ $\neq$ $t_{k_0}^0+h$.
\item \it{Case 5} : $k_0$ $\in$ $\{1,\ldots,m-1\}$, $t_0$ $=$ $t_{k_0}^0 +h$.
\end{itemize}

\noindent $\tri$ {\it Cases 1 and 2}~: these two cases imply directly from   \reff{condtermcomp1} (resp.  \reff{condtermcomp2})  that 
$M\leq 0$.

\vspace{1mm}

\noindent $\tri$ {\it Cases 3 and 4}~: we focus only on  case 4, as  case 3 involves similar (and simpler) arguments. 
We follow general viscosity solution technique based on the Ishii technique and work towards a contradiction. 
To this end, let us consider the following function~:
\beqs
\Phi_\eps(t,t',x,x',p,p') &:=& \overline{u_{k_0}}(t,x,p) - \underline{w_{k_0}}(t',x',p') - \psi_\eps(t,t',x,x',p,p'),
\enqs
with 
\beqs
\psi_\eps(t,t',x,x',p,p') &=& \frac{1}{2}\big[|t-t_0|^{2} + |p-p_0|^2 \big]  +  \frac{1}{4} |x-x_0|^{4}  \\
& & \;\; +  \frac{1}{2\eps}\big[ [t-t'|^2 + |x-x'|^2 + [p-p'|^2\big]. 
\enqs
 By the positiveness of the function $\psi_\eps$, we notice that $(t_0,x_0,p_0)$ is a strict maximizer of
$(t,x,p)$ $\rightarrow$ $\Phi_{\eps}(t,t,x,x,p,p)$. Hence, by Proposition 3.7 in \cite{craishlio92}, there exists a sequence of maximizers 
$(t_{\eps},t'_{\eps},x_{\eps},x'_{\eps},p_{\eps},p'_{\eps})$ of $\Phi_{\eps}$ such that~:
\beq
(t_\eps,t'_\eps,x_\eps,x'_\eps,p_\eps,p'_\eps) & \to &  (t_0,t_0,x_0,x_0,p_0,p_0),  \label{convAlpha}\\
\overline{u_{k_0}}(t_\eps,x_\eps,p_\eps) - \underline{w_{k_0}}(t'_\eps,x'_\eps,p'_\eps) & \to &
\overline{u_{k_0}}(t_0,x_0,p_0) - \underline{w_{k_0}}(t_0,x_0,p_0), \label{convFonct} \\
\frac{1}{\eps}\big[ |t_\eps-t'_\eps|^2 + |x_\eps-x'_\eps|^2 + |p_\eps -p'_\eps|^2] &   \to &  0 \; \mbox{ as } \;  \eps \to 0. \label{prodAlpha}
\enq
By applying Theorem 3.2 in \cite{craishlio92} to the sequence of   maximizers $(t_\eps,t'_\eps,x_\eps,x'_\eps,p_\eps,p'_\eps)$ of 
$\Phi_\eps$, we
get the existence of two symmetric matrices $A_\eps,A'_\eps$ such that~:
\beq
\label{superJet}\left( r_\eps,q_\eps,A_\eps\right) \in \overline{J}^{2,+}\overline{u_{k_0}}(t_\eps,x_\eps,p_\eps)\\
\label{infJet} \left( r'_\eps,q'_\eps,A'_\eps\right) \in \overline{J}^{2,-}\underline{w_{k_0}}(t'_\eps,x'_\eps,p'_\eps),
\enq
where
\beq
r_\eps  \; = \; \Dt{\psi_\eps}(t_\eps,t'_\eps,x_\eps,x'_\eps,p_\eps,p'_\eps) & = &   \frac{1}{\eps} (t_\eps-t'_\eps)+(t_\eps-t_0), \label{defr} \\  
r'_\eps \; = \;  - \Dtp{\psi_\eps}(t_\eps,t'_\eps,x_\eps,x'_\eps,p_\eps,p'_\eps) &=&   \frac{1}{\eps} (t_\eps-t'_\eps)  \label{defr'} \\
q_\eps \; = \;  \Dx{\psi_\eps}(t_\eps,t'_\eps,x_\eps,x'_\eps,p_\eps,p'_\eps)  &=& 
\frac{1}{\eps} ( x_\eps- x'_\eps)+ \left|x_\eps-x_0\right|^{2}\left(x_\eps-x_0\right),  \label{defq}  \\
 q'_\eps \; = \;  -  \Dxp{\psi_\eps}(t_\eps,t'_\eps,x_\eps,x'_\eps,p_\eps,p'_\eps) &=&  \frac{1}{\eps}  ( x_\eps- x'_\eps),  \label{defq'} 
\enq
and
\beq \label{ineqMat}
\left(
\begin{array}{cc}
A_\eps & 0 \\
0     & -A'_\eps
\end{array}
\right) &\leq&
 \left(
\begin{array}{cc}
\frac{3}{\eps} I_{d} -Q\left(x_\eps-x_0\right) & - \frac{3}{\eps} I_{d} \\
- \frac{3}{\eps} I_{d}     & \frac{3}{\eps} I_{d} 
\end{array}
\right),
\enq
with
\beqs
Q(x) &=&  2 x \otimes x +|x|^{2}I_{d}, 
\enqs
$I_{d}$  the identity matrix of dimension $d\times d$, and for $x$ $=$ $(x_i)_{1\leq i\leq d}$ $\in$ $\R^d$, $x \otimes x$ is the tensorial product defined by 
$x \otimes x=(x_{i}x_{j})_{i,j \in \left\{1..d\right\}^{2}}$.  Here, to alleviate notations, and since there is no derivatives with respect to  the variable $p$ in the 
PDE,  the semi-jets are defined with respect to the variables $(t,x)$, and we omitted the terms corresponding to the derivatives of 
$\psi_\eps$ with respect to $p$.   
%We also omitted the terms corresponding to $t_\eps$ and $p_\eps$ in $A_\eps$, which would only complicate the notations, 
%without introducing any new  difficulties. Indeed, one could derive the rest of proof filling the corresponding terms in the volatility matrix with zeros. 
We set $p_\eps$ $=$ $(t_i^\eps,e_i^\eps)_{1\leq i\leq k_0}$, and  $p'_\eps$ $=$ $(t_i^{'\eps},e_i^{'\eps})_{1\leq i\leq k_0}$. 
 From \reff{convAlpha}, we deduce that for $\eps$ small enough,
$t_\eps$ $\in$ $\T_{p_0}^2(k_0)$ and  $t_\eps$ $\neq$ $t_{k_0}^\eps+h$.
From  \reff{superJet}-\reff{infJet}, and the formulation of viscosity subsolution of $u_{k_0}$ to \reff{viscoQ2} and $\eta$-strict viscosity 
supersolution of $w_{k_0}$ to \reff{viscoQ2} by means of semi-jets,   
we have  for all $\eps$ small enough~:
 \beq
\min\left\{ - r_\eps - b(x_\eps) q_\eps - \frac{1}{2} {\rm tr} \left(\sigma\sigma'(x_\eps)A_\eps \right)-f(x_\eps),\right. & & \nonumber \\
\hspace{13mm} \left. \overline{u_{k_0}}(t_\eps,  x_\eps,p_\eps) - \sup_{e\in E}\overline{u_{k_0+1}}(t_\eps,x_\eps,p_\eps\cup(t_\eps,e))\right\}
&\leq& 0,  \label{sousSolJet} \\ 
 \min\left\{ - r'_\eps - b(x'_\eps) q'_\eps - \frac{1}{2} {\rm tr} \left(\sigma\sigma'(x'_\eps)A_\eps' \right)-f(x'_\eps),\right. & &\nonumber \\
 \hspace{13mm}  \left. \underline{w_{k_0}}(t'_\eps,  x'_\eps,p'_\eps)
- \sup_{e\in E}\underline{w_{k_0+1}}(t'_\eps,x'_\eps,p'_\eps\cup(t'_\eps,e))\right\}&\geq& \eta.  \label{surSolJet} 
 \enq
We then distinguish the following two possibilities in  \reff{sousSolJet}~: 
\begin{itemize}
 \item (i)  for all $\eps$ small enough, 
\beqs 
\overline{u_{k_0}}(t_\eps,  x_\eps,p_\eps) - \sup_{e\in E}\overline{u_{k_0+1}}(t_\eps,x_\eps,p_\eps\cup(t_\eps,e)) &\leq& 0. 
 \enqs
 Then, for all $\eps$ small enough, there exists $e_{\eps}$ $\in$ $E$ such that~:
 \beqs
 \overline{u_{k_0}}(t_\eps,  x_\eps,p_\eps) &\leq& \overline{u_{k_0+1}}(t_\eps,x_\eps,p_\eps\cup(t_\eps,e_\eps)) + \frac{\eta}{2}.
 \enqs
Moreover, by \reff{surSolJet}, we have
\beqs
\underline{w_{k_0}}(t'_\eps,  x'_\eps,p'_\eps) &\geq&  \underline{w_{k_0+1}}(t'_\eps,x'_\eps,p'_\eps\cup(t'_\eps,e_{\eps})) + \eta.
\enqs
Combining the two above inequalities, we deduce that for all $\eps$ small enough, 
\beqs
& & \overline{u_{k_0}}(t_\eps,  x_\eps,p_\eps)-\underline{w_{k_0}}(t'_\eps,  x'_\eps,p'_\eps)  \\
&\leq &  \overline{u_{k_0+1}}(t_\eps,x_\eps,p_\eps\cup(t_\eps,e_\eps)) - \underline{w_{k_0+1}}(t'_\eps,x'_\eps,p'_\eps\cup(t'_\eps,e_\eps)) - \frac{\eta}{2}.
\enqs
Since  $E$ is compact, there exists some $e$ $\in$ $E$ s.t. $e_\eps \to e$ up to a subsequence. From \reff{convAlpha}-\reff{convFonct}, and since 
$\overline{u_{k_0}}$, $-\underline{w_{k_0}}$ are u.s.c., we obtain by sending $\eps$ to zero~:
\beqs
& &\overline{u_{k_0}}(t_0,  x_0,p_0)-\underline{w_{k_0}}(t_0,  x_0,p_0)  \\
& \leq &  \overline{u_{k_0+1}}(t_0,x_0,p_0\cup(t_0,e)) - \underline{w_{k_0+1}}(t_0,x_0,p_0\cup(t_0,e)) - \frac{\eta}{2},
\enqs
which contradicts \reff{defM}.

\item  (ii)  for all $\eps$ small enough, 
\beqs
- r_\eps - b(x_\eps) q_\eps - \frac{1}{2} {\rm tr} \left(\sigma\sigma'(x_\eps)A_\eps \right)-f(x_\eps) &\leq& 0.  \label{diffOpPos}
\enqs
Combining with \reff{surSolJet},  we then get
\beq
\eta & \leq & r_\eps-r'_\eps + b(x_\eps) q_\eps-b(x'_\eps) q'_\eps \nonumber \\
& & \;\;\; + \;  \frac{1}{2} {\rm tr} \left(\sigma\sigma'(x_\eps)A_\eps-\sigma\sigma'(x'_\eps)A'_\eps \right) + f(x_\eps)-f(x'_\eps). \label{etacon} 
\enq
We now analyze the convergence of  the r.h.s. of \reff{etacon}  as $\eps$ goes to zero.  
First, we see from  \reff{convAlpha} and \reff{defr}-\reff{defr'} that $r_\eps-r'_\eps$ converge to  zero.  We also immediately see from the continuity of $f$ and 
\reff{convAlpha} that $f(x_\eps)-f(x'_\eps)$ converge to zero.   It is also clear from the Lipschitz property of $b$, \reff{convAlpha}, \reff{prodAlpha}, and 
\reff{defq}-\reff{defq'}  that $b(x_\eps) q_\eps-b(x'_\eps)q'_\eps$ converge to zero.  Finally, from \reff{ineqMat}, we have 
\beqs
{\rm tr} \left(\sigma\sigma'(x_\eps)A_\eps-\sigma\sigma'(x'_\eps)A'_\eps \right) &\leq&
\frac{3}{\eps} {\rm tr} \left((\sigma(x_\eps)-\sigma(x'_\eps))(\sigma(x_\eps)-\sigma(x'_\eps))'\right)\\
& & \;\;\; -  \;  {\rm tr} \left(\sigma\sigma'(x_\eps)Q(x_\eps-x_0)\right),
\enqs
and the r.h.s. of the above inequality converges to zero from the Lipschitz property of $\sigma$, \reff{convAlpha} and \reff{prodAlpha}.  
Therefore, by sending $\eps$ to zero into \reff{etacon}, we obtain the required contradiction~:  $\eta$ $\leq$ $0$. 
\end{itemize}

\noindent $\tri$  {\it Case 5}~: We keep the same notations as in the previous case.  
The crucial difference is that $\overline{u_{k_0}}$ and $\underline{w_{k_0}}$ may be sub and supersolution to different equations, depending on 
the position of $t_\eps$ (resp. $t'_\eps$) with respect to  $t_{k_0}^\eps+h$ (resp. $t_{k_0}^{'\eps}+h$). Actually, up to a subsequence for $\eps$, we 
have three subcases.  If  $t_\eps$ $\geq$ $t_{k_0}^\eps+h$ and $t_\eps'$ $\geq$ $t_{k_0}^{'\eps}+h$ for all $\eps$ small enough,  
the proof of the preceding case applies.  If  $t_\eps$ $<$ $t_{k_0}^\eps+h$,  for all $\eps$ small enough, then we have  the viscosity subsolution (resp. supersolution)  property of  $\overline{u_{k_0}}$ (resp. $\underline{w_{k_0}}$)  to  the same linear PDE~:  $-\Dt{v_k} - \Lc v_k   - f $ $=$ $0$, 
at  $(t_\eps,x_\eps,p_\eps)$ (resp.  $(t'_\eps,x'_\eps,p'_\eps)$),  and we conclude as in {\it Case 3}.    
Finally,  if $t_\eps$ $\geq$ $t_{k_0}^\eps+h$ and $t'_\eps$ $<$  $t_{k_0}^{'\eps}+h$ for all $\eps$ small enough,  then the viscosity subsolution property of 
$u_{k_0}$ to \reff{viscoQ2} at $(t_\eps,x_\eps,p_\eps)$, and the viscosity $\eta$-strict supersolution property of $w_{k_0}$ to \reff{viscoQ1} at 
$(t'_\eps,x'_\eps,p'_\eps)$ lead to~: 
\beq \label{r'geq}
- r'_\eps - b(x'_\eps) q'_\eps - \frac{1}{2} {\rm tr} \left(\sigma\sigma'(x'_\eps)A_\eps' \right)-f(x'_\eps)&\geq& \eta
\enq
and the following two possibilities~: 
\beq
- r_\eps - b(x_\eps) q_\eps - \frac{1}{2} {\rm tr} \left(\sigma\sigma'(x_\eps)A_\eps \right)-f(x_\eps)&\leq& 0,   \label{rleq} 
\enq
or
\beq
\overline{u_{k_0}}(t_\eps,  x_\eps,p_\eps) - \sup_{e\in E}\overline{u_{k_0+1}}(t_\eps,x_\eps,p_\eps\cup(t_\eps,e)) &\leq & 0.  \label{uPlusGros} 
\enq
The first possibility  \reff{r'geq}, \reff{rleq} is dealt with by the same arguments as in {\it Case 4} (ii).  The second possibility  
\reff{r'geq}, \reff{uPlusGros} does not allow to conclude directly. In fact, we  use the  additional condition \reff{condSuppl}~: 
\beq
\underline{w_{k_0}}(t_0,x_0,p_0) &\geq& \sup_{e\in E}\underline{w_{k_0+1}}(t_0,x_0,p_0\cup(t_0,e)) + \eta.
\enq
Since $\underline{w_{k_0}}$ is lower semicontinuous, this implies by \reff{convAlpha} that for all $\eps$ small enough~: 
\beqs
\underline{w_{k_0}}(t'_\eps,x'_\eps,p'_\eps) &\geq& \underline{w_{k_0}}(t_0,x_0,p_0) - \frac{\eta}{2} \\
&\geq& \sup_{e\in E}\underline{w_{k_0+1}}(t_0,x_0,p_0\cup(t_0,e)) + \frac{\eta}{2}.
\enqs
Hence, by combining  with \reff{uPlusGros}, we deduce that
\beqs
& & \overline{u_{k_0}}(t_\eps,  x_\eps,p_\eps) - \underline{w_{k_0}}(t'_\eps,x'_\eps,p'_\eps) + \frac{\eta}{2} \\
&\leq& \sup_{e\in E}\overline{u_{k_0+1}}(t_\eps,x_\eps,p_\eps\cup(t_\eps,e))- \sup_{e\in E}\underline{w_{k_0+1}}(t_0,x_0,p_0\cup(t_0,e)), 
\enqs
for all $\eps$ small enough.  From  \reff{convFonct} and Lemma \ref{lemmeVathana}, we  then obtain by sending $\eps$ to zero~:
\beqs
& & \overline{u_{k_0}}(t_0,x_0,p_0) - \underline{w_{k_0}}(t_0,x_0,p_0) + \frac{\eta}{2} \\
& \leq & \sup_{e\in E} \overline{u_{k_0+1}}(t_0,x_0,p_0\cup(t_0,e)) -  \sup_{e\in E}\underline{w_{k_0+1}}(t_0,x_0,p_0\cup(t_0,e)) \\
&\leq& \sup_{e\in E}\left\{\overline{u_{k_0+1}}(t_0,x_0,p_0\cup(t_0,e))
- \underline{w_{k_0+1}}(t_0,x_0,p_0\cup(t_0,e))\right\}. 
\enqs
This is in contradiction with \reff{defM}. 
\ep

\vspace{3mm}

Finally,  as usual,  the comparison theorem for strict supersolutions implies comparison for supersolutions.

\vspace{2mm}

\noindent {\bf Proof  of Proposition \ref{procompn}}

\noindent For any $\eta>0$, we use Lemma  \ref{existStrict} to obtain an $\eta$-strict supersolution
$w_{k}^{\eta}$ of \reff{viscoQ1}-\reff{viscoQ2}, which satisfies \reff{encadr}, so that $\underline{w_k}(t,x,p)$ $\rightarrow$
$\underline{w_k^\eta}(t,x,p)$ for all
$(t,x,p)$ $\in$ $\Dc_k^m$, as $\eta$ goes to zero. We then use Lemma \ref{comparStict} to deduce that
$\overline{u_{k}}$ $\leq$ $\underline{w_{k}}^{\eta}$ on $\Dc_k^m(n)$,  $k$ $=$ $0,\ldots,m$.
Thus, letting $\eta \to 0$, completes the proof.
\ep

\subsection{Boundary data and continuity} \label{paracont}

In this paragraph, we shall derive by induction  the boundary data \reff{vkvk-1}-\reff{vkT} in Proposition \ref{propdata}, and the continuity
of the value functions  as byproducts of viscosity properties and sequential comparison principles.

We first show relation \reff{vkT}, which follows easily from the definition of the value functions.

\begin{Lemma} \label{lem1term}
(i) For $k$  $=$ $1,\ldots,m$, $p \in \Theta_k(m)\times E^k$, $(t,x)$ $\in$ $\T_p(k)\times\R^d$,  we have~:
\beq \label{vkthetam}
v_k(t,x,p) &=& \E\Big[ \int_t^T f(X_s^{t,x,0}) ds + g(X_T^{t,x,0}) \Big].
\enq
(ii) Relation \reff{vkthetam} also holds for $k$ $=$ $0$, and for all  $(t,x)$ $\in$ $\T^m(0)\times\R^d$. In particular, $v_0(T^-,x)$ exists and is equal to $v_0(T,x)$ $=$ $g(x)$.
\end{Lemma}
{\bf Proof.}
(i) Fix $k$ $=$ $1,\ldots,m$, $p$ $=$ $(t_i,e_i)_{1\leq i\leq k}$ $\in$ $\Theta_k(m)\times E^k$, and $(t,x)$ $\in$
$\T_p(k)\times\R^d$ $=$ $[t_k,T]\times\R^d$. By definition of $\Theta_k(m)$, we have $t_i+mh$ $>$ $T$, $i$ $=$ $1,\ldots,k$.
Then, for all $\alpha$ $=$ $(\tau_i,\xi_i)_{i\geq 1}$ $\in$ $\Ac_{t,p}$, we have $\tau_i+mh$ $>$ $T$, so that from \reff{eqX}, $X_s^{t,x,p,\alpha}$ $=$
$X_s^{t,x,0}$ for $t\leq s\leq T$. We deduce immediately \reff{vkthetam} from the definition of $v_k$.

\noindent (ii) This assertion was already stated in \reff{v0f0} as a consequence of \reff{vkthetam} for $k$ $=$ $1$ and \reff{viscoQ2} for $k$ $=$ $0$.
\ep

\vspace{2mm}

The derivation of relation \reff{vkvk-1} is more delicate. We first state the following result, which is a direct consequence of the dynamic programming principle.

\begin{Lemma}
\noindent (i) For $k$ $=$ $1,\ldots,m$, and $p$ $=$ $(t_i,e_i)_{1\leq i\leq k}$ $\in$ $\Theta_k^m\times E^k$,  we have for all $x$ $\in$ $\R^d$, and $t$ $\in$ $\T_p^1(k)$ $=$ $[t_k,(t_k+h)\wedge(t_1+mh))$,
\beq
v_k(t,x,p) &=& \E\Big[ \int_t^{(t_k+h)\wedge (t_1+mh)} f(X_s^{t,x,0}) ds   + v_k(t_k+h,X_{t_k+h}^{t,x,0},p) 1_{t_k+h < t_1+mh} \label{vkt11}   \\
& & \;  + \; \Big(  c(X_{t_1+mh}^{t,x,0},e_1) +  v_{k-1}(t_1+mh,\Gamma(X_{t_1+mh}^{t,x,0},e_1),p_-) \Big) 1_{t_1+mh\leq t_k+h} \Big].  \nonumber
\enq
(ii) For $k$ $=$ $1,\ldots,m$, and $p$ $=$ $(t_i,e_i)_{1\leq i\leq k}$ $\in$ $\Theta_k^m\times E^k$, such that
$t_k+h$ $<$ $t_1+mh$, we have  for all $x$ $\in$ $\R^d$, and $t$ $\in$ $\T_p^2(k)$ $=$ $[t_k+h,t_1+mh)$,
\beq
v_k(t,x,p) & \geq & \E\Big[ \int_t^{t_1+mh} f(X_s^{t,x,0}) ds  \nonumber \\
& & \;\;\;\;\;\;\; + \; c(X_{t_1+mh}^{t,x,0},e_1) +  v_{k-1}(t_1+mh,\Gamma(X_{t_1+mh}^{t,x,0},e_1),p_-) \Big] \label{vkt121} \\
v_k(t,x,p)  & \leq &  \sup_{(\tau,\xi)\in\Ic_t} \E\Big[ \int_t^{(t_1+mh)\wedge\tau} f(X_s^{t,x,0}) ds
+ v_{k+1}(\tau,X_\tau^{t,x,0},p\cup (\tau,\xi))  1_{\tau < t_1+mh}  \nonumber \\
& &  +   \Big( c(X_{t_1+mh}^{t,x,0},e_1) +  v_{k-1}(t_1+mh,\Gamma(X_{t_1+mh}^{t,x,0},e_1),p_-) \Big) 1_{t_1+mh \leq \tau} \Big].   \label{vkt12}
\enq
\end{Lemma}
{\bf Proof.} First, we recall from  the dynamic programming principle that by making an immediate impulse control, i.e. by taking in
\reff{dynprovkpc2},  $\theta$ $=$ $t$ and $\tau$ $=$ $t$,  $\xi$ $=$ $e$ arbitrary in $E$, we have for all
$k$ $=$ $0,\ldots,m-1$, $p$ $=$ $(t_i,e_i)_{1\leq i\leq k}$ $\in$ $\Theta_k\times E^k$, $(t,x)$ $\in$ $\T_p(k)\times\R^d$ with $t$ $\geq$ $t_k+h$,
\beq \label{vksupvk+1}
v_k(t,x,p) & \geq & \sup_{e \in E} v_{k+1}(t,x,p\cup (t,e)).
\enq

\vspace{1mm}

\noindent (i) Fix $k$ $=$ $1,\ldots,m$,  $p$ $=$ $(t_i,e_i)_{1\leq i\leq k}$ $\in$ $\Theta_k^m\times E^k$,  and  $(t,x)$ $\in$ $\T_p^1(k)\times\R^d$.
We distinguish the two following cases~:

\noindent $\bullet$ {\it Case 1}~: $t_k+h$ $<$ $t_1+mh$. Then, for all $\alpha$ $\in$ $\Ac_{t,p}$, we have from \reff{eqX},
$X_s^{t,x,p,\alpha}$ $=$ $X_s^{t,x,0}$ for $t$ $\leq$ $s$ $\leq$ $t_k+h$. Hence, by applying \reff{reldynpro} with $\theta$
$=$ $t_k+h$, and noting that $\tau_i+mh$ $>$ $\theta$, $k(\theta,\alpha)$ $=$ $k$, $p(\theta,\alpha)$ $=$ $p$ for any $\alpha$ $=$ $(\tau_i,\xi_i)$ $\in$ $\Ac_{t,p}$, we obtain the required relation \reff{vkt11}, i.e.
\beqs
v_k(t,x,p) &=& \E\Big[ \int_t^{t_k+h} f(X_s^{t,x,0}) ds   + v_k(t_k+h,X_{t_k+h}^{t,x,0},p)  \Big].
\enqs

\noindent $\bullet$ {\it Case 2}~: $t_1+mh$ $\leq$ $t_k+h$. Then, for all $\alpha$ $\in$ $\Ac_{t,p}$, we have from \reff{eqX},
$X_s^{t,x,p,\alpha}$ $=$ $X_s^{t,x,0}$ for $t$ $\leq$ $s$ $<$ $t_1+mh$, and  $X_{t_1+mh}^{t,x,p,\alpha}$ $=$
$\Gamma(X_{t_1+mh}^{t,x,0},e_1)$.  Hence, by applying \reff{reldynpro} with $\theta$
$=$ $t_1+mh$, and noting that for any $\alpha$ $=$ $(\tau_i,\xi_i)$ $\in$ $\Ac_{t,p}$,  we have either
$k(\theta,\alpha)$ $=$ $k-1$, $p(\theta,\alpha)$ $=$ $p_-$ if  $\tau_{k+1}$ $>$ $t_1+mh$ (which always arises when $t_1+mh$ $<$ $t_k+h$), or
$k(\theta,\alpha)$ $=$ $k$, $p(\theta,\alpha)$ $=$ $p_-\cup (\tau_{k+1},\xi_{k+1})$  if $\tau_{k+1}$ $=$ $t_k+h$ $=$ $t_1+mh$,
we obtain
\beqs
v_k(t,x,p) &=&\sup_{\alpha\in\Ac_{t,p}}  \E\Big[ \int_t^{t_1+mh} f(X_s^{t,x,0}) ds   +  c(X_{t_1+mh}^{t,x,0},e_1) \\
& & \;\;\;\;\;\;\;  +\;   v_{k-1}(t_1+mh,\Gamma(X_{t_1+mh}^{t,x,0},e_1),p_-) 1_{\tau_{k+1} > t_1+mh}   \\
& & \;\;\;\;\;\;\;  +\; v_{k}(t_1+mh,\Gamma(X_{t_1+mh}^{t,x,0},e_1),p_- \cup (t_1+mh,\xi_{k+1}) ) 1_{\tau_{k+1} = t_1+mh=t_k+h}  \Big].
\enqs
Now, from \reff{vksupvk+1}, if $t_1+mh$ $=$ $t_k+h$, we have
$v_{k}(t_1+mh,\Gamma(X_{t_1+mh}^{t,x,0},e_1),p_- \cup (t_1+mh,\xi_{k+1}))$ $\leq$
$v_{k-1}(t_1+mh,\Gamma(X_{t_1+mh}^{t,x,0},e_1),p_-)$  for all $\xi_{k+1}$ $\Fc_{t_1+mh}$-measurable valued in $E$.   We then deduce
\beqs
v_k(t,x,p) &=&  \E\Big[ \int_t^{t_1+mh} f(X_s^{t,x,0}) ds   +  c(X_{t_1+mh}^{t,x,0},e_1)   +
v_{k-1}(t_1+mh,\Gamma(X_{t_1+mh}^{t,x,0},e_1),p_-) \Big],
\enqs
which is the required relation \reff{vkt11}.

\vspace{1mm}

\noindent (ii) Fix $k$ $=$ $1,\ldots,m$,  $p$ $=$ $(t_i,e_i)_{1\leq i\leq k}$ $\in$ $\Theta_k^m\times E^k$, s.t.
$t_k+h$ $<$ $t_1+mh$, and  $(t,x)$ $\in$ $\T_p^2(k)\times\R^d$.  Then, for all $\alpha$ $\in$ $\Ac_{t,p}$, we have from \reff{eqX},
$X_s^{t,x,p,\alpha}$ $=$ $X_s^{t,x,0}$ for $t$ $\leq$ $s$ $<$ $t_1+mh$, and  $X_{t_1+mh}^{t,x,p,\alpha}$ $=$
$\Gamma(X_{t_1+mh}^{t,x,0},e_1)$.  Let $\alpha$ $=$ $(\tau_i,\xi_i)$ be some arbitrary element in $\Ac_{t,p}$,  and set
$\tau$ $=$ $\tau_{k+1}$, $\xi$ $=$ $\xi_{k+1}$.  Observe that with $\theta$ $=$ $(t_1+mh)\wedge\tau$, we have a.s. either $k(\theta,\alpha)$ $=$ $k+1$, $p(\theta,\alpha)$ $=$ $p\cup (\tau,\xi)$ if $\tau$ $<$ $t_1+mh$ or $k(\theta,\alpha)$ $=$ $k-1$, $p(\theta,\alpha)$ $=$
$p_-$ if $\tau$ $>$ $t_1+mh$, or  $k(\theta,\alpha)$ $=$ $k$, $p(\theta,\alpha)$ $=$ $p_-\cup (\tau,\xi)$
if $\tau$ $=$ $t_1+mh$.  Hence, by applying \reff{reldynpro1} to some $\alpha$ $=$ $(\tau_i,\xi_i)$ $\in$ $\Ac_{t,p}$ s.t. $\tau_{k+1}$ $>$ $t_1+mh$ a.s. and with $\theta$ $=$ $t_1+mh$, we get the inequality \reff{vkt121}.  Furthermore, from \reff{reldynpro2}, for all $\eps$ $>$ $0$, there exists $\alpha$ $=$ $(\tau_i,\xi_i)$ $\in$ $\Ac_{t,p}$ s.t. by setting $\tau$ $=$ $\tau_{k+1}$, $\xi$ $=$ $\xi_{k+1}$, and
with $\theta$ $=$ $(t_1+mh)\wedge\tau$,
\beqs
v_k(t,x,p) -Ê\eps  & \leq &  \E\Big[ \int_t^{(t_1+mh)\wedge\tau} f(X_s^{t,x,0}) ds
+ v_{k+1}(\tau,X_\tau^{t,x,0},p\cup (\tau,\xi))  1_{\tau < t_1+mh}  \nonumber \\
& &  \; + \;  c(X_{t_1+mh}^{t,x,0},e_1) 1_{t_1+mh \leq \tau}  +   v_{k-1}(t_1+mh,\Gamma(X_{t_1+mh}^{t,x,0},e_1),p_-)
1_{t_1+mh < \tau} \\
& & \; + \;    v_{k}(t_1+mh,\Gamma(X_{t_1+mh}^{t,x,0},e_1),p_-\cup (\tau,\xi))  1_{\tau = t_1+mh} \Big].
\enqs
Now,  we have $v_{k}(t_1+mh,\Gamma(X_{t_1+mh}^{t,x,0},e_1),p_- \cup (t_1+mh,\xi))$ $\leq$
$v_{k-1}(t_1+mh,\Gamma(X_{t_1+mh}^{t,x,0},e_1),p_-)$ from \reff{vksupvk+1}.  Since $(\tau,\xi)$ $\in$ $\Ic_t$, and $\eps$ is arbitrary,
we deduce the required relation \reff{vkt12}.
\ep

\begin{Proposition} \label{procont}
For all $k$ $=$ $0,\ldots,m$,   $v_k$ is continuous on $\Dc_k^m$ and $\Dc_k(m)$. Moreover, for all $k$ $=$ $1,\ldots,m$, $p$
$=$ $(t_i,e_i)_{1\leq i\leq k}$ $\in$ $\Theta_k^m\times E^k$,  $x$ $\in$ $\R^d$,
\beqs
 v_k((t_1+mh)^-,x,p) &=& c(x,e_1) + v_{k-1}(t_1+mh,\Gamma(x,e_1),p_-).
\enqs
\end{Proposition}
{\bf Proof.}
{\bf 1.} We first easily see from Lemma \ref{lem1term}, continuity and growth condition of $f$, $g$, and dominated convergence theorem, that $v_k$ is continuous on $\Dc_k(m)$, for all $k$ $=$ $0,\ldots,m$.

\vspace{1mm}

\noindent {\bf 2.}   We shall prove by forward induction on $n$ $=$ $m+1,\ldots,N$
 that   {\bf (Hk)(n)}, $k$ $=$ $1,\ldots,m$,  and {\bf (H0)(n)} hold, where

\vspace{2mm}

{\bf (Hk)(n)}  \hspace{2mm}  $v_k$ is continuous on $\Dc_k^m(n)$, and for all  $p$ $=$ $(t_i,e_i)_{1\leq i\leq k}$
$\in$ $\Theta_k^m(n)\times E^k$,

\hspace{19mm} $v_k((t_1+mh)^-,x,p)$ $=$ $c(x,e_1) + v_{k-1}(t_1+mh,\Gamma(x,e_1),p_-)$,  $x$ $\in$ $\R^d$.

\vspace{1mm}

{\bf (H0)(n)}  \hspace{2mm} $v_0$ is continuous on $\T^n(0)\times\R^d$.

\vspace{2mm}

\noindent $\tri$ {\bf Initialization~: $n$ $=$ $m+1$.}  Let us prove that {\bf (Hk)(m+1)}, $k$ $=$ $1,\ldots,m$, and {\bf (H0)(m+1)} are satisfied.

\vspace{1mm}

\noindent $\bullet$  Take some $k$ $=$ $1,\ldots,m$, and fix some arbitrary $x$ $\in$ $\R^d$ and $p$ $=$ $(t_i,e_i)_{1\leq i\leq k}$
$\in$ $\Theta_k^m(m+1)\times\E^k$.  Notice that $p_-$ $\in$ $\Theta_{k-1}(m)\times E^{k-1}$ so that
$v_{k-1}(.,.,p_-)$ is continuous on
$\T_{p_-}(k-1)\times\R^d$  from  part  {\bf 1.} above.  Here, to alleviate notations, we used the convention that $\T_{p_-}(k-1)$ $=$ $\T^m(0)$ if $k-1$ $=$ $0$.  We distinguish two cases~:

\vspace{1mm}

\noindent $\star$  {\it Case 1.} $\T_p^2(k)$ $=$ $\emptyset$, i.e. $t_1+mh$ $\leq$ $t_k+h$ so that  $\T_p(k)$ $=$ $\T_p^1(k)$ $=$ $[t_k,t_1+mh)$. From \reff{vkt11}, we then have
for all $t$ $\in$ $\T_p(k)$~:
\beqs
v_k(t,x,p) &=& \E\Big[ \int_t^{t_1+mh} f(X_s^{t,x,0}) ds     +
c(X_{t_1+mh}^{t,x,0},e_1) +  v_{k-1}(t_1+mh,\Gamma(X_{t_1+mh}^{t,x,0},e_1),p_-)  \Big].
\enqs
By continuity of $v_{k-1}(t_1+mh,.,p_-)$, $\Gamma(.,e_1)$, $c(.,e_1)$,  growth condition on $f$, $c$, $\Gamma$ and $v_{k-1}$, we deduce with the dominated convergence theorem that $v_k((t_1+mh)^-,x,p)$ exists and
\beqs
 v_k((t_1+mh)^-,x,p) &=& c(x,e_1) + v_{k-1}(t_1+mh,\Gamma(x,e_1),p_-).
\enqs

\vspace{1mm}

\noindent $\star$  {\it Case 2.} $\T_p^2(k)$ $=$ $[t_k+h,t_1+mh)$ $\neq$ $\emptyset$, i.e. $t_1+mh$ $>$ $t_k+h$ (this implies in particular that
$k$ $<$ $m$ and $m$ $>$ $1$).  From \reff{vkt121}-\reff{vkt12}, we first  prove that
\beq
&  &  \overline{v_k}(t_1+mh,x,p) \nonumber \\
&\leq & \max\big[  c(x,e_1) + v_{k-1}(t_1+mh,x,p_-)  ,   \sup_{e \in E}   \overline{v_{k+1}}(t_1+mh,x,p\cup (t_1+mh,e))   \big].  \label{vksup}
\enq
Indeed, consider some sequence $(t_\eps,x_\eps,p_\eps)_{\eps>0}$ $\in$ $\Dc_k^m$ converging to $(t_1+mh,x,p)$  and such that
$\lim_{\eps\rightarrow 0}$ $v_k(t_\eps,x_\eps,p_\eps)$ $=$ $\overline{v_k}(t_1+mh,x,p)$. For any $\eps$ $>$ $0$, one can find, by \reff{vkt12},
some $(\hat\tau_\eps,\hat\xi_\eps)$ $\in$ $\Ic_{t_\eps}$ s.t.
\beqs
v_k(t_\eps,x_\eps,p_\eps) & \leq &  \E\Big[ \int_{t_\eps}^{(t_1^\eps+mh)\wedge\hat\tau_\eps} f(X_s^{t_\eps,x_\eps,0}) ds
+ v_{k+1}(\hat\tau_\eps,X_{\hat\tau_\eps}^{t_\eps,x_\eps,0},p_\eps\cup (\hat\tau_\eps,\hat\xi_\eps))
1_{\hat\tau_\eps < t_1^\eps+mh}  \nonumber \\
& &  + \;  \Big( c(X_{t_1^\eps+mh}^{t_\eps,x_\eps,0},e_1^\eps) +
v_{k-1}(t_1^\eps+mh,\Gamma(X_{t_1^\eps+mh}^{t_\eps,x_\eps,0},e_1^\eps),p_{\eps -}) \Big) 1_{t_1^\eps+mh \leq \hat\tau_\eps} \Big] \; + \eps,
\enqs
where we denote $p_\eps$ $=$ $(t_i^\eps,e_i^\eps)_{1\leq i\leq k}$ and $p_{\eps -}$ $=$ $(t_i^\eps,e_i^\eps)_{2\leq i\leq k}$. By setting
\beqs
G_\eps &=&
c(X_{t_1^\eps+mh}^{t_\eps,x_\eps,0},e_1^\eps)
+  v_{k-1}(t_1^\eps+mh,\Gamma(X_{t_1^\eps+mh}^{t_\eps,x_\eps,0},e_1^\eps),p_{\eps -}),
\enqs
we  rewrite the above inequality as
\beq
v_k(t_\eps,x_\eps,p_\eps) & \leq & \E\Big[ \int_{t_\eps}^{(t_1^\eps+mh)\wedge\hat\tau_\eps}
f(X_s^{t_\eps,x_\eps,0}) ds  \; + \;  G_\eps \nonumber \\
& & \;\;\; + \; \Big(v_{k+1}(\hat\tau_\eps,X_{\hat\tau_\eps}^{t_\eps,x_\eps,0},p_\eps\cup (\hat\tau_\eps,\hat\xi_\eps)) - G_\eps\Big)
1_{\hat\tau_\eps < t_1^\eps+mh} \Big] \; + \eps.  \label{vkeps}
\enq
Since $p_-$ $\in$ $\Theta_{k-1}(m)\times E^{k-1}$, we have $p_{\eps -}$ $\in$ $\Theta_{k-1}(m)\times E^{k-1}$ for $\eps$ small enough.
Hence, by continuity of $v_{k-1}$ on $\Dc_{k-1}(m)$ (from part {\bf 1.}), continuity of $\Gamma$ and $c$, and path-continuity of the flow
$X_s^{t,x,0}$, we have
\beq \label{limGeps}
\lim_{\eps\rightarrow 0} G_\eps &=& G \; := \;  c(x,e_1) + v_{k-1}(t_1+mh,\Gamma(x,e_1),p_-) \;\;\;\;\; a.s.
\enq
Moreover, by compactness of $E$, the sequence $(\hat\xi_\eps)_\eps$ converges, up to a subsequence, to some $\xi$ valued in $E$.
We deduce that
\beq
& & \limsup_{\eps\rightarrow 0}  \Big(v_{k+1}(\hat\tau_\eps,X_{\hat\tau_\eps}^{t_\eps,x_\eps,0},p_\eps\cup (\hat\tau_\eps,\hat\xi_\eps)) - G_\eps\Big)
1_{\hat\tau_\eps < t_1^\eps+mh} \nonumber \\
& \leq & \Big(\overline{v_{k+1}}(t_1+mh,x,p\cup (t_1+mh,\xi)) - G\Big) \limsup_{\eps\rightarrow 0}  1_{\hat\tau_\eps < t_1^\eps+mh}
\nonumber \\
& \leq &  \Big(\sup_{e\in E} \overline{v_{k+1}}(t_1+mh,x,p\cup (t_1+mh,e)) - G\Big) \limsup_{\eps\rightarrow 0}
1_{\hat\tau_\eps < t_1^\eps+mh} \;\;\; a.s. \label{limvk+1eps}
\enq
From the linear growth condition on $f$, $c$, $\Gamma$, $v_{k-1}$, $v_{k+1}$, and estimate \reff{estimX2},  we may use dominated convergence theorem and send $\eps$ to zero in \reff{vkeps} to obtain with \reff{limGeps}-\reff{limvk+1eps}~:
\beqs
& & \overline{v_k}(t_1+mh,x,p) \\
& \leq & \E\Big[ G + \Big(\sup_{e\in E} \overline{v_{k+1}}(t_1+mh,x,p\cup (t_1+mh,e)) - G\Big) \limsup_{\eps\rightarrow 0}
1_{\hat\tau_\eps < t_1^\eps+mh} \Big] \\
& \leq & \max\Big[ G, \sup_{e\in E} \overline{v_{k+1}}(t_1+mh,x,p\cup (t_1+mh,e))\Big],
\enqs
which is the required inequality \reff{vksup}.

We next show that
\beq \label{inegvk+1}
\sup_{e \in E}   \overline{v_{k+1}}(t_1+mh,x,p\cup (t_1+mh,e)) & \leq & c(x,e_1) + v_{k-1}(t_1+mh,x,p_-).
\enq
Indeed, for any arbitrary $e$ $\in$ $E$, consider some sequence $(t_\eps,x_\eps,p_\eps,e_\eps)_{\eps>0}$ $\in$
$\Dc_k^m\times E$ converging to $(t_1+mh,x,p,e)$  and such that
$\lim_{\eps\rightarrow 0}$ $v_{k+1}(t_\eps,x_\eps,p_\eps\cup (t_\eps,e_\eps))$ $=$
$\overline{v_{k+1}}(t_1+mh,x,p\cup (t_1+mh,e))$.  For $\eps$ small enough, $t_\eps+h$ $\geq$ $t_1^\eps+mh$, and so from the DPP \reff{vkt11}, we have~:
\beq
v_{k+1}(t_\eps,x_\eps,p_\eps\cup (t_\eps,e_\eps)) &=& \E\Big[ \int_{t_\eps}^{t_1^\eps+mh} f(X_s^{t_\eps,x_\eps,0}) ds +
c(X_{t_1^\eps+mh}^{t_\eps,x_\eps,0},e_1^\eps) \nonumber \\
& & \;\;\;\;\;  + \; v_k(t_1^\eps+mh,\Gamma(X_{t_1^\eps+mh}^{t_\eps,x_\eps,0},e_1^\eps),p_{\eps -} \cup (t_\eps,e_\eps)) \Big].
\label{vk+1inter2}
\enq
Since $p_-$ $\in$ $\Theta_{k-1}(m)\times E^{k-1}$, we have $p_{\eps -}$ $\in$ $\Theta_{k-1}(m)\times E^{k-1}$ for $\eps$ small enough.
Hence, by continuity of $v_{k}$ on $\Dc_{k}(m)$, continuity and growth linear  condition of $f$, $\Gamma$ and $c$, and path-continuity of the flow
$X_s^{t,x,0}$, we send $\eps$ to zero in \reff{vk+1inter2} and get  by the dominated convergence theorem
\beq
& &  \overline{v_{k+1}}(t_1+mh,x,p\cup (t_1+mh,e))  \nonumber \\&=&
c(x,e_1) + v_k (t_1+mh,\Gamma(x,e_1),p_- \cup (t_1+mh,e)). \label{vkvk-1e}
\enq
Moreover, from \reff{vksupvk+1},  we have $v_{k-1}(t_1+mh,\Gamma(x,e_1),p_-)$ $\geq$ $v_k(t_1+mh,\Gamma(x,e_1),p_-\cup (t_1+mh,e))$ for all $e$ $\in$ $E$. Plugging into \reff{vkvk-1e}, this proves \reff{inegvk+1}.

Finally, we easily see from \reff{vkt121} that
\beq \label{vkinf}
\underline{v_k}(t_1+mh,x,p)  & \geq   & c(x,e_1) + v_{k-1}(t_1+mh,x,p_-).
\enq
Indeed, consider some sequence $(t_\eps,x_\eps,p_\eps)_{\eps>0}$ $\in$ $\Dc_k^m$ converging to $(t_1+mh,x,p)$  and such that
$\lim_{\eps\rightarrow 0}$ $v_k(t_\eps,x_\eps,p_\eps)$ $=$ $\underline{v_k}(t_1+mh,x,p)$. From \reff{vkt121}, we have in particular
\beqs
v_k(t_\eps,x_\eps,p_\eps) & \geq & \E\Big[ \int_{t_\eps}^{t_1^\eps+mh} f(X_s^{t_\eps,x_\eps,0}) ds \\
& & \;\;\;\;\;  + \;
c(X_{t_1^\eps+mh}^{t_\eps,x_\eps,0},e_1^\eps) +
v_{k-1}(t_1^\eps+mh,\Gamma(X_{t_1^\eps+mh}^{t_\eps,x_\eps,0},e_1^\eps),p_{\eps -}) \Big].
\enqs
By continuity and linear growth condition of $v_{k-1}$,  $\Gamma$, $c$, $f$, and estimate \reff{estimX2},  we get  \reff{vkinf} by  the dominated convergence theorem, and  sending $\eps$ to zero in the above inequality. 

Hence, the inequalities \reff{vksup}-\reff{inegvk+1}-\reff{vkinf} prove that  $v_k((t_1+mh)^-,x,p)$ exists and is equal to~:
\beq
v_k((t_1+mh)^-,x,p) &=&   \overline{v_k}(t_1+mh,x,p) \; = \;  \underline{v_k}(t_1+mh,x,p) \label{vklim} \\
&=& c(x,e_1) + v_{k-1}(t_1+mh,\Gamma(x,e_1),p_-).  \nonumber
\enq
We have then proved that \reff{vklim} holds for all $k$ $=$ $1,\ldots,m$, $p$ $=$ $(t_i,e_i)_{1\leq i\leq k}$ $\in$
$\Theta_k^m(m+1)\times E^k$, and $x$ $\in$ $\R^d$.

\vspace{1mm}

\noindent $\bullet$  We know from Proposition \ref{provisco} that the family of value functions $v_k$, $k$ $=$ $0,\ldots,m$, is
a viscosity solution to \reff{viscoQ1}-\reff{viscoQ2}, in particular at step $n$ $=$ $m+1$. We also recall from Lemma \ref{lem1term} that $\overline{v_0}(T,x)$ $=$ $\underline{v_0}(T,x)$ $=$ $g(x)$. Together with \reff{vklim}, and the comparison principle at step $n$ $=$ $m+1$ in Proposition \ref{procompn}, this proves $\overline{v_k}$ $\leq$ $\underline{v_k}$ on $\Dc_k^m(n)$ for $n$ $=$ $m+1$.
This implies the continuity of $v_k$ on $\Dc_k(m+1)$, $k$ $=$ $0,\ldots,m$, and so {\bf (Hk)(m+1)}, $k$ $=$ $1,\ldots,m$, and {\bf (H0)(m+1)} are stated.

\vspace{2mm}

\noindent $\tri$ {\bf Step $n$  $\rightarrow$  $n+1$ for $n$ $\in$ $\{m+1,\ldots,N-1\}$}. We suppose that
{\bf (Hk)(n)}, $k$ $=$ $1,\ldots,m$, and {\bf (H0)(n)} hold true.

\vspace{1mm}

\noindent  Take some $k$ $=$ $1,\ldots,m$, and fix some arbitrary $x$ $\in$ $\R^d$ and $p$ $=$ $(t_i,e_i)_{1\leq i\leq k}$
$\in$ $(\Theta_k^m(n+1)\times\E^k$.  Notice that $p_-$ $\in$ $\Theta_{k-1}(n)\times E^{k-1}$.  By same arguments as in step $n$ $=$ $m+1$, using here,
instead of part {\bf 1.},  continuity of $v_{k-1}$ on $\Dc_{k-1}^m(n)$ by the induction hypothesis of step $n$, we prove that
\beqs
v_k((t_1+mh)^-,x,p) &=&   \overline{v_k}(t_1+mh,x,p) \; = \;  \underline{v_k}(t_1+mh,x,p)  \\
&=& c(x,e_1) + v_{k-1}(t_1+mh,\Gamma(x,e_1),p_-).  \nonumber
\enqs
We also have $\overline{v_0}(T,x)$ $=$ $\underline{v_0}(T,x)$ $=$ $g(x)$. Therefore, from the viscosity property of $v_k$, $k$ $=$ $0,\ldots,m$, to \reff{viscoQ1}-\reff{viscoQ2} at step $n+1$, and the comparison principle in at step $n+1$
in Proposition \ref{procompn}, we obtain $\overline{v_k}$ $\leq$ $\underline{v_k}$ on $\Dc_k^m(n+1)$, which implies the continuity of
$v_k$ on  $\Dc_k^m(n+1)$, $k$ $=$ $0,\ldots,m$. Therefore, {\bf (Hk)(n+1)}, $k$ $=$ $1,\ldots,m$, and {\bf (H0)(n+1)} are proved.

\vspace{2mm}

\noindent $\tri$ The proof is completed at step $N$ by recalling  that $\Theta_k(N)$ $=$ $\Theta_k$, $\Dc_k^m(N)$ $=$ $\Dc_k^m$, $k$ $=$ $0,\ldots,m$.
\ep

\subsection{Proof of Theorem \ref{thmmain}}

In view of the  results proved in paragraphs \ref{paravisco} and \ref{paracont}, it remains to prove the uniqueness  result of Theorem \ref{thmmain}.
Let us then consider another family $w_k$, $k$ $=$ $0,\ldots,m$  of viscosity solutions to \reff{viscoQ1}-\reff{viscoQ2}, satisfying growth condition
\reff{growthvk}, and boundary data \reff{vkvk-1}-\reff{vkT}~:  for $k$ $=$ $1,\ldots,m$, $p$ $=$ $(t_i,e_i)_{1\leq i\leq k}$ $\in$ $\Theta_k^m\times E^k$,
$x$ $\in$ $\R^d$,
\beq \label{wkwk-1}
w_k((t_1+mh)^-,x,p) &=& c(x,e_1) + w_{k-1}(t_1+mh,\Gamma(x,e_1),p_-).
\enq
and
\beq \label{wkT}
w_k(t,x,p) &=&   \E\Big[ \int_t^T f(X_s^{t,x,0}) ds + g(X_T^{t,x,0}) \Big], \;\;\; (t,x,p) \in \Dc_k(m).
\enq
We shall prove by forward induction on $n$ $=$ $m,\ldots,N$ that $v_k$ $=$ $w_k$ on $\Dc_k(n)$.

\vspace{2mm}

\noindent $\tri$ {\bf Initialization~: $n$ $=$ $m$}.  Relations \reff{vkT}, \reff{wkT} and  \reff{v0f0}  show that $v_k$ $=$ $w_k$
on $\Dc_k(m)$, $k$ $=$ $0,\ldots,m$.

\vspace{1mm}

\noindent $\tri$ {\bf Step $n$ $\rightarrow$ $n+1$}.  Suppose that $v_k$ $=$ $w_k$ on $\Dc_k(n)$, $k$ $=$ $0,\ldots,m$.
For any $k$ $\geq$ $1$, $p$ $=$ $(t_i,e_i)_{1\leq i\leq k}$ $\in$ $\Theta_k^m(n+1)\times E^k$, we notice that $p_-$ $\in$
$\Theta_{k-1}(n)\times E^{k-1}$. Hence $v_{k-1}(t_1+mh,\Gamma(x,e_1),p_-)$ $=$ $w_{k-1}(t_1+mh,\Gamma(x,e_1),p_-)$, $x$ $\in$ $\R^d$,
and so from \reff{vkvk-1}, \reff{wkwk-1}, we have
\beqs
v_k((t_1+mh)^-,x,p) &=& w_k((t_1+mh)^-,x,p).
\enqs
We already know that $v_0(T^-,x)$ $=$  $w_0(T^-,x)$ ($=$ $g(x)$).   Therefore, from the comparison principle at step $n+1$ in
Proposition \ref{procompn}, we deduce that $u_k$ $=$ $w_k$ on $\Dc_k^m(n+1)$, and so on $\Dc_k(n+1)$, $k$ $=$ $0,\ldots,m$.  Finally, the proof is completed since $\Dc_k(N)$ $=$ $\Dc_k$.

\end{document}